\documentclass[preprint, 11pt]{elsarticle}
\usepackage{lineno}
\pdfoutput=1

\usepackage[utf8]{inputenc}
\usepackage{amsmath}
\usepackage{amssymb}
\usepackage{mathrsfs}
\usepackage{color}
\usepackage{xcolor}
\usepackage{graphicx}
\usepackage[percent]{overpic}
\usepackage{url}
\usepackage{caption}
\usepackage{subcaption}
\usepackage{enumerate}
\usepackage{overpic}
\usepackage{amsthm}
\usepackage{mathtools}
\usepackage{moreverb}
\usepackage{multirow}
\usepackage{multicol}
\usepackage[pdftex,colorlinks,bookmarksopen,bookmarksnumbered,citecolor=red,urlcolor=red]{hyperref}
\usepackage{mwe}
\usepackage{graphbox}
\usepackage{textcomp}
\usepackage{tabularx}

\def\bu{{\bm u}}

\def\x{{\bm x}}

\newcommand{\cF}{\mathcal F}
\newcommand{\cI}{\mathcal I}
\newcommand{\cL}{\mathcal L}
\def\0{\boldsymbol{0}}

\def\dt{\partial_t}

\def\cl {\nonumber \\}

\newtheorem{rem}{Remark}[section]

\newcommand{\bm}[1]{\mbox{\boldmath{$#1$}}}
\def\div{\nabla\cdot}

\newcommand{\anna}[2][cyan]{{\textcolor{#1}{#2}}}

\usepackage{verbatim}
\usepackage{graphicx}
\usepackage{epstopdf}
\journal{Journal of Computational and Applied Mathematics}

\begin{document}

%% ELSEVIER %%
\begin{frontmatter}

%%%% arxiv %%%
%\date{}
%\title{A novel Large Eddy Simulation model for the Quasi-Geostrophic Equations in a Finite Volume setting}
%\author[1]{Michele Girfoglio\thanks{mgirfogl@sissa.it}}
%\author[2]{Annalisa Quaini\thanks{quaini@math.uh.edu}}
%\author[1]{Gianluigi Rozza\thanks{grozza@sissa.it}}
%\affil[1]{SISSA, International School for Advanced Studies, Mathematics Area, mathLab, via Bonomea, Trieste 265 34136, %Italy}
%\affil[2]{Department of Mathematics, University of Houston, Houston TX 77204, USA}
%\maketitle

\title{
Nonlinear filtering stabilizations 
for the quasi-geostrophic equations
}
%\tnotetext[mytitlenote]{Fully documented templates are available in the elsarticle package on \href{http://www.ctan.org/tex-archive/macros/latex/contrib/elsarticle}{CTAN}.}

%% Group authors per affiliation:
%\author{Elsevier\fnref{myfootnote}}
%\address{Radarweg 29, Amsterdam}
%\fntext[myfootnote]{Since 1880.}

%% or include affiliations in footnotes:
%\author[mymainaddress,mysecondaryaddress]{Elsevier Inc}
%\ead[url]{www.elsevier.com}

%\author[mysecondaryaddress]{Global Customer Service\corref{mycorrespondingauthor}}
%\cortext[mycorrespondingauthor]{Corresponding author}
%\ead{support@elsevier.com}

%\address[mymainaddress]{1600 John F Kennedy Boulevard, Philadelphia}
%\address[mysecondaryaddress]{360 Park Avenue South, New York}

\author[Clemson]{Lander Besabe}
\ead{lbesabe@clemson.edu}

\author[Houston]{Sachin Kumar}
\ead{skumar43@Cougarnet.uh.edu}

\author[TUM]{Annalisa Quaini}
\ead{annalisa.quaini@tum.de}

% \author[SISSA]{Gianluigi Rozza}
% \ead{grozza@sissa.it}
% %\ead[url]{www.elsevier.com}

\address[Clemson]{School of Mathematical and Statistical Sciences, Clemson University, Clemson, SC 29634-0975, USA}

\address[Houston]{Department of Mathematics, University of Houston, Houston TX 77204, USA}

\address[TUM]{Chair of Computational Methods in Civil and Environmental Engineering, Technical University of Munich,
Arcisstra{\ss}e 21, 80333 Munich, Germany}

\begin{abstract}
Numerical simulations of ocean flows typically require fine computational meshes to resolve the Munk scale, leading to high computational costs. Filtering-based large eddy simulation (LES) provides a way to relax the mesh size requirement by modeling the effects of the unresolved scales. For the implementation of this strategy, we 
propose a three-step algorithm called Evolve-Filter-Relax (EFR) that requires (i) the solution of a QGE problem, (ii) a nonlinear Helmholtz filter for the
potential vorticity field leveraging an indicator function, and (iii) a final relaxation step.
We show that the EFR algorithm can be interpreted as a splitting scheme for a perturbed QGE problem with additional dissipation and provide a practical choice for the relaxation parameter. For comparison, we also investigate a nonlinear Bardina regularization of the QGE.
Numerical results on a classical benchmark show that both the EFR approach and the nonlinear Bardina regularization 
significantly improve the accuracy and stability of coarse mesh simulations with no LES model.
Additionally, the EFR method
with a deconvolution-based indicator function 
delivered
the best balance between accuracy, stability, and computational efficiency in a test
case involving a more realistic geometry (Mediterranean Sea).
\end{abstract}


\end{frontmatter}

%\linenumbers

\centerline{\emph{Dedicated to the memory of Alessandro Veneziani}}

\section{Introduction}\label{sec:intro}

%\url{https://agupubs.onlinelibrary.wiley.com/doi/full/10.1029/2018JC014636}

%\anna{TO DO}

Ocean circulation plays a critical role in regulating the Earth's climate through heat distribution \cite{Trenberth2001} and freshwater transport \cite{Talley2008}. Further, the ocean influences the planet's climate variability across a wide range of spatial and temporal scales \cite{Wunsch2004, Griffies2000}. Since observing the ocean directly is a difficult task, the numerical simulation of ocean flows has been essential in studying its dynamics, understanding its interaction with other components of the climate, and improving climate change predictions and responses \cite{Fox-Kemper2019, Marshall1997, Marshall1997b}.
However, the simulation of large-scale ocean circulation remains a challenging task due to the size of the problem (surface area of order $1e+06$ Km$^2$ and depth of up to $11$ Km$^2$) and the typically long time intervals of interest. Ocean flows are also characterized by the evolution 
of flow structures with a broad range of spatial scales (as typical for high Reynolds numbers), requiring very fine computational mesh to resolve the smaller scales. This poses
an additional challenge at the computational level.
Thus, much effort has been devoted to reducing
the computational cost of ocean flow simulations. 
See, e.g., \cite{Greatbatch2000, Foster2013, QGE-review}.

In this work, we will focus on one of the simplest models for mid-latitude wind driven ocean dynamics: the quasi-geostrophic equations (QGE). The name refers to the fact that in ocean motion the Coriolis force and pressure gradient forces are almost balanced, but inertia also has an effect.
The Rossby number $Ro$ is a dimensionless number which characterizes the strength of inertia compared to the strength of the Coriolis force. The QGE hold in the limit of small $Ro$, so that inertial forces are orders of magnitude smaller than the Coriolis and pressure forces.
Although the QGE are a simplified model, they retain the features that make ocean flow simulations computationally challenging. 
In fact, when 
the Munk scale (a length that depends on the ratio of $Ro$ and the Reynolds number) is small, the simulations require very fine meshes because the mesh size has to be smaller than the Munk scale.

To reduce the computational cost, simulations are often performed with coarse meshes and an additional model, called subgrid-scale model (SGS), to account for the unresolved scales. A common approach is to introduce a constant (in time and space) eddy viscosity that is substantially larger than the kinetic viscosity of water \cite{Bryan1963,Gates1968,Holland1975, Berloff1999, Tanaka2010}. %The required spatial resolution is then closely related to the Munk scale, a quantity that depends on both $Re$ and $Ro$. When the computational mesh does not resolve the Munk scale, additional stabilization is generally required to obtain an accurate and stable solution. 
This could be considered a rough 
Large Eddy Simulation (LES) approach. 
A more sophisticated approach would 
introduce an artificial viscosity that evolves in 
space and time in response to the
flow features and the underlying mesh. 
See, e.g., \cite{sagaut2006large, BIL05, layton2012approximate}.

Several LES strategies have been developed and investigated for the QGE. One of the earliest works \cite{Nadiga2001} used dispersive-dissipative eddy  parametrization for the barotropic QGE. In \cite{Holm2003}, an $\alpha$-regularization was introduced for the same model. 
A related regularization approach is the BV-Bardina model, which was introduced and analyzed in \cite{Monteiro2015}. 
Approximate deconvolution LES was applied to the single-layer QGE in \cite{San2011} and subsequently extended to the two-layer case in \cite{San2012}. A multi-scale coarse grid projection technique was developed in \cite{San_IJMCE2013} to improve the accuracy and stability of the QGE solution on a given coarse mesh. In \cite{Maulik2016}, Smagorinsky and Leith closure models were used together with an efficient trapezoidal filter. 
% I removed the next two sentences because it felt like I was adding too many details for this particular work. We mentioned this robustness result in 4.1.
%The artificial viscosity in the Leith model is proportional to the gradient of the vorticity field, but in the Smagorinsky model, it is related to the local strain rate field. It is concluded that the Leith model adds less eddy viscosity and is less robust to changes in mesh sizes, as opposed to the Smagorinsky model which is robust to changes in mesh sizes. 
A modified SGS model for the QGE was developed in \cite{Babu2026} by combining the Leith model with an approximate deconvolution-based nonlinear gradient model to balance stabilization and representation of the backscattering of energy. 
More recently, differential filtering has been developed as a stabilization strategy for the under-resolved QGE. In \cite{Girfoglio_JCAM2023}, linear and nonlinear differential filters were applied to the potential vorticity, where the nonlinear filter is based on the gradient of the potential vorticity. This technique has also been considered for reduced order models of the one-layer \cite{Girfoglio2023} and 
two-layer \cite{Besabe2025} QGE.

An alternative filtering-based LES is the Evolve-Filter-Relax (EFR) framework. At each time step, the governing equations are first evolved on a coarse mesh, then the resulting solution is filtered, and a relaxation step combines the filtered and unfiltered solutions. 
Introduced in \cite{layton_JMFM}, EFR algorithms
have been widely applied to the incompressible Navier-Stokes equations \cite{BQV,Bowers2012,abigail_CMAME,Ervin2012,layton_JMFM,Girfoglio2019,LAYTON20113183,Olshanskii2013}. It was shown in \cite{Bowers2012,abigail_CMAME,Ervin2012,layton_JMFM,LAYTON20113183} that numerical results obtained with nonlinear differential filters are more precise in localizing where eddy viscosity is needed and are overall more accurate than results obtained with plain Smagorinsky-type models.
More recently, the EFR framework has been applied to the weakly compressible Euler equations \cite{Clinco_2023}. Its application to the QGE remains unexplored and, 
to the best of our knowledge, 
it is considered for the first time in this paper. 

This paper shows that the EFR algorithm can be interpreted as a splitting scheme for a perturbed QGE with additional dissipation and we derive a practical choice for the relaxation parameter $\chi$ based on the heuristic approach from \cite{BQV}. 
Linear and nonlinear filters, including Leith-, Smagorinsky-, and deconvolution-based indicator functions, 
are investigated
and their effects on accuracy and stability of the system are assessed. %We consider two challenging flow regimes using the classical double-gyre wind forcing benchmark \cite{Greatbatch2000}, and a more challenging physical test involving a wind-driven flow in the Mediterranean sea. 
In a series of numerical tests, 
the LES solutions are compared with high-resolution reference solutions in terms of their time-averaged potential vorticity and stream function, kinetic energy evolution, and kinetic energy spectra. %Finally, we investigate the sensitivity of EFR to the filtering radius parameter. 
For comparison, we also investigate the performance of linear and nonlinear variants of the BV-Bardina model.

%consider two filtering-based LES approaches for coarse simulations of the QGE: the BV-Bardina model and the Evolve-Filter-Relax (EFR) framework. The BV-Bardina model is closely related to the BV-$\alpha$ model \cite{Holm2003} and is obtained from a simplification of the Bardina model under the $\beta$-plane approximation \cite{Khouider2008, Monteiro2015}. BV-Bardina regularizes the nonlinear convective term through the filtered potential vorticity and can capture backscaterring of energy \cite{Horiuti1989}. In contrast, EFR introduces artifical dissipation through a modular three-step procedure which has been successfully applied to convection-dominated flows \cite{BQV, Girfoglio2019,Clinco_2023}. 
%We investigate linear filtering together with several nonlinear indicator functions, including Leith-, Smagorinsky-, and deconvolution-based functions, and assess their effects on the accuracy and stability of the QG system. 

The rest of the paper is organized as follows. Sec.~\ref{sec:QGE} describes the QGE model. In Sec.~\ref{sec:les}, we discuss 
BV-$\alpha$ and BV-Bardina models and present 
a nonlinear extension of the BV-Bardina model. Sec.~\ref{sec:qge-efr-desc} introduces the QGE-EFR model. In Sec.~\ref{sec:num_res},  we present several numerical results to test the performance of the LES models considered in this work. Conclusions are drawn in Sec.~\ref{sec:conclusion}.

\section{Problem definition}
\label{sec:QGE}

The QGE describe the evolution of a rotating fluid layer of constant depth $H$, temperature, and density $\rho$ under the assumptions of near hydrostatic and geostrophic balance.
%First, we assume that the ocean flow is 
%near hydrostatic and geostrophic balance.
An additional assumption in this model is
the so called $\beta$-plane approximation: 
the Coriolis frequency $\tilde f$ is linearized as
$\tilde f = \tilde f_0 + \beta \tilde y$, where $\tilde f_0$ is the local rotation rate at the center of the basin $\tilde y = 0$ (vertical coordinate) and $\beta$ is the gradient of the Coriolis frequency. 
The value of $\tilde f_0$ depends on the rotational speed of the earth and the latitude at $\tilde y = 0$. 
The $\beta$-plane approximation is equivalent to approximating the Earth (a sphere) with a tangent plane at $\tilde y = 0$. Let $\Omega=[\tilde x_0, \tilde x_f] \times [-L/2, L/2]$ be a computational domain on this plane, with coordinates $\tilde x-\tilde y$. We will consider the formulation of the 
QGE in terms of potential vorticity $\tilde q$ and stream function $\tilde \psi$.

It is common to state this model in non-dimensional form. For this purpose, let $V$ be the Sverdrup velocity \cite{San2012} 
\begin{equation*}
    V = \frac{\pi \tau_0}{\rho H \beta L},
\end{equation*}
where $\tau_0$ is the maximum amplitude of the double-gyre wind stress. 
We introduce
the non-dimensional variables: 
\begin{equation}
    x = \frac{\tilde{x}}{L}, \quad y = \frac{\tilde{y}}{L}, \quad t = \frac{\tilde{t}}{L/V}, \quad q = \frac{\tilde{q}}{\beta L}, \quad \psi = \frac{\tilde{\psi}}{VL}.
    \label{eq:change}
\end{equation}
In the non-dimensionalization of the QGE, two well-known non-dimensional numbers appear: the Reynolds number $Re$ and the Rossby number $Ro$
\begin{equation}\label{eq:ro_re}
Re = \frac{VL}{\nu}, \quad Ro = \frac{V}{\beta L^2},
\end{equation}
where $\nu$ is the eddy viscosity. 
Calling the above $Re$ Reynolds number is an abuse
of terminology since the actual 
Reynolds number is defined with the molecular viscosity of the fluid instead of the eddy viscosity.

The non-dimensional QGE read:
Find $q(x, y, t)$ and $\psi(x,y,t)$ such that
\begin{align}
    &\frac{\partial q}{\partial t} + \nabla \cdot ( (\nabla \times \bm{\psi})q) - \frac{1}{Re}\Delta q = F \quad &&\mbox{ in }\Omega \times (t_0,T), \label{eq:qge1}\\
    &q = -Ro \Delta \psi + y \quad &&\mbox{ in }\Omega \times (t_0,T), \label{eq:qge2}
\end{align}
where $\bm{\psi} = (0, 0, \psi)$ and $F$ is external forcing due to the wind. We recall that the non-dimensional potential vorticity $q$ is defined as $q = Ro~\omega + y$, where
$\omega$ is the nondimensional vorticity.

To close problem \eqref{eq:qge1}-\eqref{eq:qge2}, proper boundary conditions and initial data should be provided. 
Following \cite{Nadiga2001,Holm2003,Monteiro2015,Monteiro2014,San2012}, we enforce
\begin{align}
\psi &= 0 \quad \mbox{ on }\partial \Omega \times (t_0,T), \label{eq:BV5_comp} \\
q &= y \quad \mbox{ on }\partial \Omega \times (t_0,T), \label{eq:BV5_comp2} \\
q(x,y,t_0) &= y \quad \mbox{ in } \Omega. \label{eq:BV5_comp3}
\end{align}

For a chosen eddy viscosity, 
it is typical in the literature
\cite{Bryan1963, Gates1968, Holland1975, Berloff1999, BERLOFF_KAMENKOVICH_PEDLOSKY_2009, Tanaka2010} to consider a Direct Numerical Simulation (DNS) of the QGE a simulation with a mesh size $h$ smaller than the Munk scale \cite{Munk1950}
\begin{equation}\label{eq:re_ro}
    \delta_M = L\left(\frac{Ro}{Re}\right)^{1/3}.
\end{equation}
Even when the Reynolds number is computed with the artificial viscosity as in \eqref{eq:re_ro}, 
resolving the Munk scale often leads to high computational cost. 
As mentioned in Sec.~\ref{sec:intro},
one way to contain the computational cost is by using an LES approach: one coarsens the mesh, i.e., $h > \delta_M$, and adds an SGS model to account for the effects of the unresolved scales onto the resolved one. 
This is typically the result of a space filtering operation. 

In the next section, we discuss two LES models presented in the literature: the BV-$\alpha$ model and the BV-Bardina model. For the latter, we introduce a nonlinear variant. 

\section{BV-$\alpha$ and BV-Bardina models} \label{sec:les}

A well studied LES model for the QGE is the so called BV-\emph{$\alpha$} model (see, e.g.,  \cite{Nadiga2001, Holm2003, Monteiro2015, Monteiro2014}), whose name
owes to the fact that the 
QGE are also known as barotropic vorticity (BV) equations. The starting point for this model
is a type of turbulence closure called Lagrangian-averaged Navier–Stokes alpha (LANS-$\alpha$)
equations, which
%. The alpha
%terms in this model
%do \emph{not} introduce additional dissipative terms beyond those
%in the Navier–Stokes equations. Instead, 
alters the nonlinearity in the 
Navier-Stokes equations
so that excitations
at length scales smaller than $\alpha$ (given
length scale) are nonlinearly swept by motions at the
larger length scales.
%In this model, the transport velocity is filtered relative to the
%transported velocity as done in Leray-$\alpha$ model \cite{Leray1934}. 
Moreover, an additional nonlinear
alpha term %beyond the Leray component 
is included to restore the Kelvin circulation
theorem \cite{Chen1998}.
By taking the curl of the momentum equation in the
LANS-$\alpha$ model, reducing it to two planar dimensions with the $\beta$-plane approximation, one obtains the so-called 
BV-$\alpha$ model \cite{Holm2003}. 

The non-dimensional BV-$\alpha$ model reads:
Find $q(x, y, t)$, $\psi(x,y,t)$, and $\overline{q}(x,y,t)$ such that
\begin{align}%\label{eq:leray}
\frac{\partial q}{\partial t} + \div \left(\left(\nabla \times \bm{\psi}\right) q \right) - \dfrac{1}{Re} \Delta q &= F \quad \mbox{ in }\Omega \times (t_0,T), \label{eq:BV1}\\
-\alpha^2\Delta \overline{q} +\overline{q} &= q  \quad ~~ {\rm in}~\Omega \times
(t_0,T), \label{eq:BV2} \\
-Ro \Delta \psi + y &= \overline{q} \quad ~ \mbox{ in }\Omega \times (t_0,T), \label{eq:BV3}
\end{align}
where $\overline{q}$ is the \emph{filtered vorticity} and $\alpha$ is the above-mentioned prescribed length that can be interpreted as a filtering radius. Problem \eqref{eq:BV1}-\eqref{eq:BV3} is supplemented with boundary conditions \eqref{eq:BV5_comp}-\eqref{eq:BV5_comp3} plus
the additional boundary condition 
\begin{align}
\overline{q} &= y \quad \quad \quad \mbox{ on }\partial \Omega \times (t_0,T). \label{eq:BC_qbar}
\end{align}
In \cite{Girfoglio_JCAM2023}, we introduced a following nonlinear version of the 
BV-$\alpha$ model.

%\subsection{The BV-Bardina model} \label{sec:bvb-desc}

The BV-Bardina model is %a LES model 
obtained by simplifying the Bardina regularization to the QGE under the $\beta$-plane approximation \cite{Khouider2008, Monteiro2015}. Unlike simple eddy-viscosity closures (e.g., Smagorinsky), which introduce dissipation through an artificial viscosity, the Bardina model regularizes the nonlinear advection term by replacing the transported vorticity with its filtered counterpart. The Bardina model has been shown to capture backscattering of energy (i.e., the transfer of energy from small to large scales) \cite{Horiuti1989}, which is of particular relevance in geophysical flows \cite{Juricke2020,Sroka_Guimond_2021}.

The nondimensional BV-Bardina problem reads: Find $q(x, y, t)$, $\psi(x,y,t)$, and $\overline{q}(x,y,t)$ such that
\begin{align}%\label{eq:leray}
\frac{\partial q}{\partial t} + \div \left(\left(\nabla \times \bm{\psi}\right) \overline{q} \right) - \dfrac{1}{Re} \Delta q &= F \quad \mbox{ in }\Omega \times (t_0,T), \label{eq:BV1_bard} %\\
%-\alpha^2\Delta \overline{q} +\overline{q} &= q  \quad ~~ {\rm in}~\Omega \times
%(t_0,T), \label{eq:BV2_bard} \\
%-Ro \Delta \psi + y &= \overline{q} \quad ~ \mbox{ in }\Omega \times (t_0,T). \label{eq:BV2_bard}
\end{align}
and \eqref{eq:BV2}-\eqref{eq:BV3} hold.
The BV-Bardina model is closely related to the BV-$\alpha$ model \eqref{eq:BV1}-\eqref{eq:BV3},
but it introduces stronger regularization. 
In \cite{Monteiro2015}, it was shown that 
the BV-Bardina model gives a
good coarse mesh approximation to the highly resolved direct numerical simulation of the QGE, 
and it is more accurate than the BV-$\alpha$ model
for a well-known benchmark test. 

In \cite{Girfoglio_JCAM2023}, we showed that 
linear filter \eqref{eq:BV2}, while mathematically
and computationally convenient, has limited 
effectivity when very coarse meshes are considered. 
We found that increased accuracy can be achieved by replacing \eqref{eq:BV2} with nonlinear filter
\begin{align}%\label{eq:leray}
-\alpha^2 \div ( a(q,\psi) \nabla\overline{q}) +\overline{q} &= q  \quad {\rm in}~\Omega \times
(t_0,T), \label{eq:BV2_bard_NL}
\end{align}
where $0<a(\cdot, \cdot)\leq 1$ is an indicator function used to learn where and how much artificial viscosity is needed. In particular, the indicator function takes values close to zero where the flow field does not need regularization and close to one where the flow field does need regularization.
While we write $a$ as a function of both $q$ and $\psi$ to be general, it could be a function of either $q$ or $\psi$.
See Sec.~\ref{sec:ind} for the indicator functions considered in this work.

The nonlinear extension of the BV-Bardina model: Find $q(x, y, t)$, $\psi(x,y,t)$, and $\overline{q}(x,y,t)$ such that \eqref{eq:BV3},\eqref{eq:BV1_bard},\eqref{eq:BV2_bard_NL} hold. 
% \begin{align}%\label{eq:leray}
% \frac{\partial q}{\partial t} + \div \left(\left(\nabla \times \bm{\psi}\right) \overline{q} \right) - \dfrac{1}{Re} \Delta q &= F \quad \mbox{ in }\Omega \times (t_0,T), \label{eq:BV1_bard_NL}\\
% -\alpha^2 \div ( a(q) \nabla\overline{q}) +\overline{q} &= q  \quad {\rm in}~\Omega \times
% (t_0,T), \label{eq:BV2_bard_NL} \\
% -Ro \Delta \psi + y &= \overline{q} \quad ~ \mbox{ in }\Omega \times (t_0,T). \label{eq:BV2_bard}
% \end{align}

We discretize in time and linearize the nonlinear BV-Bardina model using a segregated algorithm. To present it, let $\Delta t \in \mathbb{R}$, $t^n = t_0 + n \Delta t$, with $n = 0, ..., N_T$ and $T = t_0 + N_T \Delta t$. We denote by $f^n$ the approximation of a generic quantity $f$ at the time $t^n$. The algorithm reads:
given $q^n$ and $\bm{\psi}^n$, at $t^{n+1}$ perform the following steps:
 \begin{itemize}
  \item[i)] Find the filtered vorticity $\overline{q}^{n+1}$ such that
  \begin{align}
 -\alpha^2\div \left(a^{n} \nabla\overline{q}^{n+1}\right) + \overline{q}^{n+1} = q^{n}, \label{eq:BV1_Bard_Ntd} 
 \end{align}
 where $a^{n}  = a(q^{n}, \psi^n)$.
  \item[ii)] Find the stream function $\psi^{n+1}$ such that
 \begin{align}
 -\text{Ro} \Delta \psi^{n+1} + y = \overline{q}^{n+1}. \label{eq:BV2_Bard_Ntd}
 \end{align}
 \item[iii)] Find the vorticity $q^{n+1}$ such that
 \begin{align}
\dfrac{1}{\Delta t} q^{n+1} + \div \left(\left(\nabla \times \bm{\psi}^{n+1} \right) \overline{q}^{n+1} \right) - \dfrac{1}{\text{Re}} \Delta q^{n+1} = F + \dfrac{1}{\Delta t} q^{n}. \label{eq:BV3_Bard_Ntd} %\quad \mbox{ in }\Omega \times (t_0,T), \label{eq:BV1_comp11} \\
     \end{align} 
 \end{itemize}

 \begin{rem}\label{rem:subit}
The presence of the filtered vorticity $\overline{q}$ in the nonlinear convective term of \eqref{eq:BV1_bard} introduces a tighter coupling between $\psi$ and $q$
than in the BV-$\alpha$ models. As a consequence, the loosely coupled algorithm \eqref{eq:BV1_Bard_Ntd}-\eqref{eq:BV3_Bard_Ntd}
might be unstable for highly convection-dominated flows and certain choices of discretization parameters. In such cases, one can introduce the following fixed-point subiterations (with index $k$): at time $t^{n+1}$, iteration $k+1$, given $q^{n+1,k}$ and $\bm{\psi}^{n+1,k}$, perform the following steps
\begin{itemize}
  \item[i)] Find the filtered vorticity $\overline{q}^{n+1,k+1}$ such that
 \begin{align}
 -\alpha^2\div \left(a^{n+1,k} \nabla\overline{q}^{n+1,k+1}\right) + \overline{q}^{n+1,k+1} = q^{n+1,k}, \label{eq:BV1_Bard_Ntd_it} 
 \end{align}
 where $a^{n+1,k}  = a(q^{n+1,k}, \psi^{n+1,k})$.
  \item[ii)] Find the stream function $\psi^{n+1,k+1}$ such that
 \begin{align}
 -\text{Ro} \Delta \psi^{n+1,k+1} + y = \overline{q}^{n+1,k+1}. \label{eq:BV2_Bard_Ntd_it}
 \end{align}
 \item[iii)] Find the vorticity $q^{n+1,k+1}$ such that
 \begin{align}
\dfrac{1}{\Delta t} q^{n+1,k+1} + \div \left(\left(\nabla \times \bm{\psi}^{n+1,k+1} \right) \overline{q}^{n+1,k+1} \right) - \dfrac{1}{\text{Re}} \Delta q^{n+1,k+1} = F + \dfrac{1}{\Delta t} q^{n+1,k}. \label{eq:BV3_Bard_Ntd_it} %\quad \mbox{ in }\Omega \times (t_0,T), \label{eq:BV1_comp11} \\
     \end{align} 
      \item[iv)] Check stopping criterion
      \begin{equation}\label{eq:stop}
    \frac{\|q^{n+1, k+1}-q^{n+1, k}\|_{\infty}}{\|q^{n+1, k+1}\|_{\infty}}\leq \epsilon,
\end{equation}
where $\epsilon$ is a user-defined tolerance. If satisfied, set $\overline{q}^{n+1} = \overline{q}^{n+1,k+1}$, $\psi^{n+1} = \psi^{n+1,k+1}$, and $q^{n+1} = q^{n+1,k+1}$. If not, go back to step i). 
 \end{itemize}
Clearly, the improved stability properties of algorithm
\eqref{eq:BV1_Bard_Ntd_it}-\eqref{eq:BV3_Bard_Ntd_it} comes at the cost of additional computational time. 

The $L^\infty$-norm in \eqref{eq:stop} is employed because some indicator functions take much larger values within boundary layers than in the interior of the domain and thus a global norm, such as the $L^2$-norm, might be misleading. See Sec.~\ref{sec:bv-bardina} for more details on this.
 \end{rem}

\section{The Evolve-Filter-Relax algorithm for the QGE} \label{sec:qge-efr-desc}

Unlike the BV-$\alpha$ and BV-Bardina models, which incorporates filtering into the governing equations, the Evolve-Filter-Relax (EFR) framework introduces regularization as an additional module to the problem. The method consists of three sequential steps: an \textit{evolve} step which solves the QGE without regularization, a \textit{filter} step, and a \textit{relax} step which mitigates possible over-smoothing. We will call this method QGE-EFR.

The QGE-EFR algorithm reads as follows: given $q^n$ and $\bm{\psi}^n$, at $t^{n+1}$ perform the following steps
  \begin{itemize}
 \item[i)] \emph{Evolve}: Find intermediate potential vorticity $p^{n+1}$ and stream function $\psi^{n+1}$ such that
 \begin{align}
\dfrac{1}{\Delta t} p^{n+1} + \div \left(\left(\nabla \times \bm{\psi}^n \right) p^{n+1} \right) - \dfrac{1}{\text{Re}} \Delta p^{n+1} = F + \dfrac{1}{\Delta t} q^{n}, \label{eq:BV1_td} \\
 -\text{Ro} \Delta \psi^{n+1} + y = {p}^{n+1}. \label{eq:BV2_td}
    \end{align}
    To improve computational efficiency, eq.~\eqref{eq:BV1_td} has been decoupled
    from eq.~\eqref{eq:BV2_td} by making the convective velocity explicit.
%where we have replaced $\bm{\psi}^{n+1}$ in \eqref{eq:BV1_NL} by $\bm{\psi}^n$, i.e.~a linear extrapolation. 
 \item[ii)] \emph{Filter}: Find the filtered vorticity $\overline{p}^{n+1}$ such that
  \begin{align}
 -\alpha^2\div \left(a^{n+1} \nabla\overline{p}^{n+1}\right) + \overline{p}^{n+1} = p^{n+1}, \label{eq:BV_EFR2} 
 \end{align}
 where $a^{n+1}  = a(p^{n+1}, \psi^{n+1})$. 
 \item[iii)] \emph{Relax}: For $0 \leq \chi \leq 1$, set the end-of-step potential vorticity
 \begin{align}
     q^{n+1} = (1 - \chi) p^{n+1} + \chi \overline{p}^{n+1}. \label{eq:relax}
 \end{align}
 \end{itemize}

The filtering step \eqref{eq:BV_EFR2} may be rewritten as
   \begin{align}
   \overline{p}^{n+1} = \mathcal{F} (p^{n+1}), \quad \cF = (\cI + \cL)^{-1}, \quad \cL=-\nabla \cdot (\alpha^2 a^{n+1} \nabla) \label{eq:filter} 
 \end{align}
where $\cI$ is the identity operator.
By shifting the index $n+1$ to $n$ in \eqref{eq:relax}-\eqref{eq:filter} and plugging them into \eqref{eq:BV1_td}, we obtain:
\begin{align}
\dfrac{1}{\Delta t} ( p^{n+1} - p^{n}) + \div \left(\left(\nabla \times \bm{\psi}^n \right) p^{n+1} \right) - \dfrac{1}{\text{Re}} \Delta p^{n+1} + \frac{\chi}{\Delta t} (\cI -  \cF) p^{n}= F, \label{eq:BV_EFR1} 
\end{align}
Eq.~\eqref{eq:BV_EFR1} gives us an implicit discretization of  \eqref{eq:qge1} with the Backward Euler scheme and an additional,
explicitly treated, dissipation term (the last term at the left-hand side). 
Let us assume that $\chi = \chi_0 \Delta t$,
where $\chi_0$ is a time-independent constant.
Then, \eqref{eq:BV_EFR1} can be seen as a time-stepping scheme for problem:
\begin{equation}\label{eq:model_LES}
\frac{\partial q}{\partial t} + \nabla \cdot ( (\nabla \times \bm{\psi})q) - \frac{1}{Re}\Delta q + {\chi_0}(\cI -  \cF) q= F.
\end{equation}
Thus, the EFR algorithm \eqref{eq:BV1_td}-\eqref{eq:relax} can be interpreted as a splitting scheme for problem \eqref{eq:qge2},\eqref{eq:model_LES}, which can be read as problem \eqref{eq:qge1}-\eqref{eq:qge2} with an additional term designed to be dissipative, as is typical of LES. Since it is equivalent to an LES model, we will also refer to algorithm \eqref{eq:BV1_td}-\eqref{eq:relax} as the QGE-EFR model. 

Next, let us consider the space discretization of 
steps \eqref{eq:BV1_td}-\eqref{eq:relax} and show the perturbation terms that the EFR algorithm adds
to the QGE.
We rewrite \eqref{eq:BV1_td} as
 \begin{align}
\dfrac{p^{n+1}_h - q^{n}_h}{\Delta t}  + \div \left(\left(\nabla \times \bm{\psi}^n_h \right) p^{n+1}_h \right) - \dfrac{1}{\text{Re}} \Delta p^{n+1}_h = F_h, \label{eq:BV1_td_sd}
\end{align}
where the subscript $h$ denotes a space-discrete quantity. 
Let us also write the space-discrete version of \eqref{eq:relax}
 \begin{align}
     q^{n+1}_h = (1 - \chi) p^{n+1}_h + \chi \overline{p}^{n+1}_h. \label{eq:relax_sd}
 \end{align}
 and \eqref{eq:BV_EFR2} divided by $\Delta t$
 \begin{align}
 \dfrac{\overline{p}^{n+1}_h - p^{n+1}_h}{\Delta t}
 -\div \left( \overline{\nu}_h\nabla\overline{p}^{n+1}_h\right)  = 0, \label{eq:BV_EFR2_sd} 
 \end{align}
 where
 \begin{align}
    \overline{\nu}_h = \frac{\alpha^2}{\Delta t}a^{n+1}_h. \label{eq:nubar_h}
\end{align}

We multiply \eqref{eq:BV_EFR2_sd} by $\chi$, add it to \eqref{eq:BV1_td_sd}, and make use of \eqref{eq:relax_sd} to obtain:
\begin{align}
\dfrac{q^{n+1}_h - q^{n}_h}{\Delta t}  + \div \left(\left(\nabla \times \bm{\psi}^n_h \right) p^{n+1}_h \right) - \dfrac{1}{\text{Re}} \Delta p^{n+1}_h - \chi \div \left( \overline{\nu}_h\nabla\overline{p}^{n+1}_h\right)= F_h, \label{eq:EFR_comb}
\end{align}
Using again \eqref{eq:relax_sd} in \eqref{eq:EFR_comb} and the space discrete version of \eqref{eq:BV2_td}, we get:
\begin{align}
&\dfrac{q^{n+1}_h - q^{n}_h}{\Delta t}  + \div \left(\left(\nabla \times \bm{\psi}^n_h \right) q^{n+1}_h \right) - \dfrac{1}{\text{Re}} \Delta q^{n+1}_h - \chi \div \left( \overline{\nu}_h\nabla\overline{p}^{n+1}_h\right) \cl
& \quad
+\chi \div \left(\left(\nabla \times \bm{\psi}^n_h \right) (p^{n+1}_h - \overline{p}^{n+1}_h) \right) - \chi \dfrac{1}{\text{Re}} \Delta (p^{n+1}_h - \overline{p}^{n+1}_h)
= F_h, \label{eq:pertQGE1} \\
& -\text{Ro} \Delta \psi^{n+1}_h + y_h - \chi (p^{n+1}_h - \overline{p}^{n+1}_h) = {q}^{n+1}. \label{eq:pertQGE2}
\end{align}
The above problem is a consistent perturbation of the original QGE 
problem because all the extra terms are multiplied by $\chi$, which vanishes as $\Delta t$ tends to zero. Indeed, recall that 
\eqref{eq:BV_EFR1} suggests taking 
$\chi = \mathcal{O}(\Delta t)$. 
Two of the three perturbation terms in \eqref{eq:pertQGE1} are clearly dissipative (first and third term multiplied by $\chi$), while the other (i.e., the second term
multiplied by $\chi$) has a convective nature. 
Further, we note that two perturbation terms
in \eqref{eq:pertQGE1} and the only extra term in \eqref{eq:pertQGE2} contain $p^{n+1}_h - \overline{p}^{n+1}_h$. From the discrete version of \eqref{eq:BV_EFR2}, it is easy to see that
 \begin{align*}
p^{n+1}_h -\overline{p}^{n+1}_h = 
 - \alpha^2 \div \left( a^{n+1}_h \nabla\overline{p}^{n+1}_h\right),
 \end{align*}
Since it is typical to set $\alpha = \mathcal{O}(h)$, $\overline{p}^{n+1}_h \rightarrow p^{n+1}_h$  as $h \rightarrow 0$.
Thus, as $\Delta t$ and $h$ tend to 0, all the perturbations in 
\eqref{eq:pertQGE1}-\eqref{eq:pertQGE2}
vanish %as follows: \anna{Do we want to include $h$ from the derivatives?}
% \begin{align}
%     \chi \div \left( \overline{\nu}_h\nabla\overline{p}^{n+1}_h\right) &\sim \alpha^2, \label{eq:perturb1} \\
%     \chi \div \left(\left(\nabla \times \bm{\psi}^n_h \right) (p^{n+1}_h - \overline{p}^{n+1}_h) \right) &\sim \alpha^2 \Delta t, \label{eq:perturb2} \\
%     \chi \dfrac{1}{\text{Re}} \Delta (p^{n+1}_h - \overline{p}^{n+1}_h) &\sim \frac{\alpha^2 \Delta t}{Re}, \label{eq:perturb3} \\
%     \chi (p^{n+1}_h - \overline{p}^{n+1}_h) &\sim \alpha^2 \Delta t, \label{eq:perturb4}
% \end{align}
and problem \eqref{eq:qge1}-\eqref{eq:qge2} is recovered. Of the three terms multiplied by $\chi$ in \eqref{eq:pertQGE1}, the first is the dominant one as the other two tend to zero faster.

\subsection{The indicator function}\label{sec:ind}

The success of filter \eqref{eq:BV2_bard_NL} depends significantly on the indicator function $a(\cdot)$, which determines where and how much filtering is needed.
In this work, we consider four indicator functions. The first is the linear indicator function given by 
\begin{equation*}
    a_{lin} = 1.
\end{equation*}
With this choice, the same amount of artificial diffusion is introduced everywhere in the computational domain and \eqref{eq:BV2_bard_NL} becomes \eqref{eq:BV2}. In contrast, the remaining three indicator functions
allow the filtering to act selectively, with values close to zero in regions that require little to no filtering and values close to one in regions that need the most filtering.

The Leith-type indicator function, given by
\begin{equation*}
    a_{Leith}(q) = \dfrac{|\nabla q|}{\text{max}\left(1, ||\nabla q||_\infty\right)},
\end{equation*}
selects regions with high potential vorticity gradient. This choice is closely related to the Leith model (hence the name), where the artificial dissipation introduced in the system is proportional to the gradient of the vorticity \cite{Maulik2016}. The Smagorinsky-type indicator:
\begin{equation*}
    a_{S}(\bu) = \dfrac{\| \nabla \bu \|}{\text{max}(1,|\nabla \bu\|_{\infty})}, \quad \text{where} \quad \bu = \left(\frac{\partial \psi}{\partial y}, - \frac{\partial \psi}{\partial x}\right),
\end{equation*}
applies more filtering in regions of large velocity gradient. Note that although we write $a_S$ as a function of $\bu$ for simplicity, it is ultimately a function of $\psi$, which justifies the notation introduced in \eqref{eq:BV2_bard_NL}.

Lastly, for the QGE-EFR, we also consider a deconvolution-based indicator function, which was shown to be particularly accurate for realistic incompressible flows \cite{BQV,Girfoglio2019}
and atmospheric flows \cite{Clinco_2023}.
This is given by
\begin{equation*}
    a_D(\bu) = \frac{\left|  \bu - \tilde{\bu} \right|}{\text{max}(1,\|\bu - \tilde{\bu} \|_{\infty})},
\end{equation*}
where $\tilde{\bu}$ is obtained from the linear Helmholtz filter
\begin{equation*}
    \tilde{\bu} - \alpha^2\Delta \tilde{\bu} = \bu  \quad {\rm in}~\Omega \times
(t_0,T),
\end{equation*}
with boundary condition $\tilde{\bu} = \bm{0}$ on $\partial \Omega \times (t_0,T)$.

\subsection{The relaxation parameter}

As shown in Sec.~\ref{sec:qge-efr-desc}, the perturbation terms introduced by the EFR algorithm 
are all multiplied by the relaxation parameter $\chi$. Thus, the amount of dissipation introduced by this method is multiplied by $\chi$. Following the heuristic approach in \cite{BQV}, we choose $\chi$ by requiring the effective dissipation on a coarse mesh with mesh size $h$ match the viscous contribution on a mesh which resolves the Munk scale $\delta_M$. %Let the subscript $h$ denote the space-discrete solution where $h$ is the mesh size, and define
%\begin{equation*}
%    \bar{\nu}_h = \frac{\alpha^2}{\Delta t}a_{h}^{n+1}.
%\end{equation*}

We use $\nabla_\zeta$ to denote the discrete gradient operator defined  on a mesh with size $\zeta$. Comparing the dominant dissipative terms introduced in the QGE-EFR algorithm:
\begin{align}
    - \nabla_h \cdot \left( \frac{1}{Re}  \nabla_h q^{n+1}_h + \chi \overline{\nu}_h  \nabla_h \overline{p}^{n+1}_h \right),\label{eq:viscous_h}
\end{align}
with the physical viscous term on a mesh of size $\delta_M$:
\begin{align}
    - \nabla_{\delta_M} \cdot \left( \frac{1}{Re} \nabla_{\delta_M} q^{n+1}_h \right), \label{eq:viscous_delta}
\end{align}
and approximating $\nabla_\zeta \sim 1/\zeta$ gives
\begin{align}
    \chi \simeq 
    \frac{1}{|| \overline{\nu}_h ||_\infty Re} \left( \frac{h^2}{\delta_M^2} - 1 \right) =
    \frac{1}{\alpha^2 || a_h^{n+1} ||_\infty Re} \left( \frac{h^2}{\delta_M^2} - 1 \right) \Delta t. \label{eq:chi}  
\end{align}
Since the quantity $\|a_h^{n+1}\|_{\infty}$ must be evaluated at every time step, we adopt the simplified expression
\begin{align*}
    \chi \simeq 
    \frac{1}{Re} \left( \frac{1}{\delta_M^2} - \frac{1}{h^2} \right) \Delta t.
\end{align*}
by taking $\alpha = h$ and knowing that $\|a_h^{n+1}\|_{\infty} \leq 1$. This formula does not depend on time and requires only a priori knowledge of the physical parameters ($Re$ and $\delta_M$) and the discretization parameters ($h$ and $\Delta t$), with the downside that it is only effective when $|| a_h^{n+1} ||_\infty$ is close to 1. For the numerical experiments in Sec.~\ref{sec:num_res}, this simplified formula is used to determine the values of $\chi$.

\section{Numerical Results} \label{sec:num_res}
In this section, we test the accuracy, efficiency, and robustness of the models presented in Sec.~\ref{sec:les}. 
For their implementation, we used open-source package Geophysical and Environmental Applications (GEA) \cite{GEA, GEAproceeding}, which is built on top of the C++ finite volume library OpenFOAM\textsuperscript{\textregistered} \cite{Weller1998}. See \cite{Girfoglio_AIP2023, Girfoglio_JCAM2023, Girfoglio2021a} for more details on GEA.

In Sec.~\ref{sec:bv-bardina} and \ref{sec:qge-efr}, we assess the performance of BV-Bardina and QGE-EFR models using the classical double-gyre wind forcing experiment, following \cite{Greatbatch2000,Holm2003,Monteiro2014,Monteiro2015}. The computational domain is a rectangular basin $\Omega = [0, 1]\times [-1, 1]$ with wind forcing $F = \sin(\pi y)$. For this test, we set $t_0 = 0$, $T = 100$, and $\Delta t = 2.5e-5$. We consider two cases:
\begin{itemize}
    \item[-] Experiment 1: $Re = 1000$ and $Ro = 0.008$ 
    \item[-] Experiment 2: $Re = 1200$ and $Ro = 0.0096$ 
\end{itemize}
We note that experiment 1 has been considered in \cite{Girfoglio_JCAM2023} to assess the performance of the linear and nonlinear BV-$\alpha$ models. Both experiments feature a nondimensional Munk scale of 
\begin{equation*}
    \delta_M/L = 0.02.
\end{equation*}
As we will see later in the section, experiment 2 is a more challenging case, and thus, more interesting, due to the larger Reynolds number. For the DNS, we employ a mesh size of $h = 1/256$, which is almost an order of magnitude smaller the Munk scale.
Note that the use of DNS here is an abuse of terminology since the Reynolds number is computed with the eddy viscosity instead of the actual viscosity. 
%\lander{LB -- I am not mentioning the fact that we take $\alpha = O(h)$ here just in case we want to include this information in the discussions preceding this section.} 
%\anna{Sounds good.}
%For each test case, simulations are performed on two coarse meshes, $h = 1/16$ and $h = 1/32$, and two filter radii, $\alpha = h$ and $\alpha = 2h$.

\begin{figure}[htb!]
    \centering
    \begin{subfigure}[h]{0.16\textwidth}
        \centering
         \includegraphics[width=\textwidth]{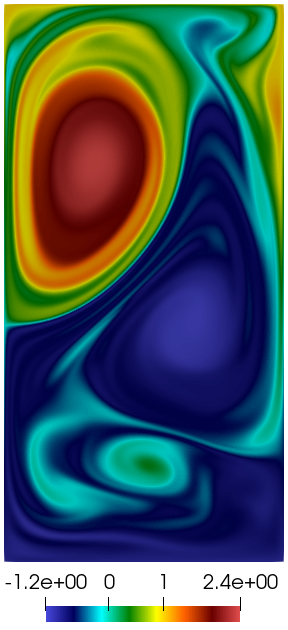}
         \caption{$q, t = 30$}
    \end{subfigure}
    \begin{subfigure}[h]{0.16\textwidth}
        \centering
         \includegraphics[width=\textwidth]{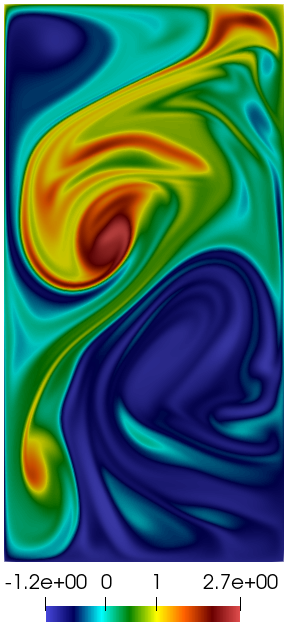}
         \caption{$q, t = 31$}
    \end{subfigure}
    \begin{subfigure}[h]{0.16\textwidth}
        \centering
         \includegraphics[width=\textwidth]{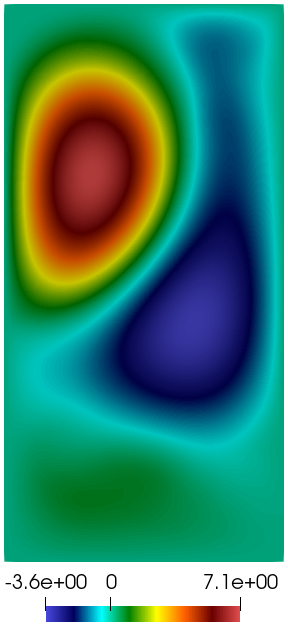}
         \caption{$\psi, t = 30$}
    \end{subfigure}
    \begin{subfigure}[h]{0.16\textwidth}
        \centering
         \includegraphics[width=\textwidth]{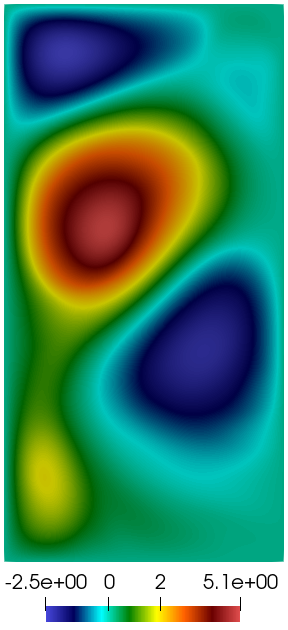}
         \caption{$\psi, t = 31$}
    \end{subfigure}
    \caption{Experiment 1: panels (a) and (b) show the potential vorticity at $t = 30$ and $t = 31$
    and panels (c) and (d) show 
    the stream functions at the same time instants.}
    \label{fig:instantaneous_fields}
\end{figure}

Fig.~\ref{fig:instantaneous_fields} shows snapshots of the potential vorticity and stream function at time instants $t = 30$ and $t = 31$ for experiment 1, demonstrating the highly convective nature of the flow. 
Hence, comparison of instantaneous snapshots of the variables between the DNS and the LES models over long periods of time is challenging and may not adequately reflect the overall performance of the LES models. 
Since it is well-known that the QGE solution reaches a statistically steady state after an initial transient phase, we focus on time-averaged quantities, $\tilde{q}$ and $\tilde{\psi}$, computed over the time interval $[20, 100]$. 
When converted to dimensional time, 
this interval corresponds to 9 years for experiment 1 and 7.5 years for experiment 2. %\lander{Sachin, will you check this information for me? I used the fact that in dimesional time, $\Delta t$ is around 3.1 minutes and we perform the average over 3.2M time steps.} \sachin{SK: According to my calculations, $\Delta t$ is approximately 3.3 minutes that means [80,100] is equivalent to approximately 20 years.}
%\anna{Is that true for both experiments?} \lander{Sachin, can you verify this for me?}
%\sachin{It is approximately 9 years for Case 1 when Re = 1000, and approximately 7.5 years for Case 2 when Re = 1200. It is approx 20 years when Re = 450.}

To account for the time-dependent behavior of the QGE system, we also look at the evolution of the kinetic energy (KE):
\begin{equation}\label{eq:ke}
    E(t) = \frac{1}{2}\int_{\Omega}\left(\left(\frac{\partial\psi}{\partial y}\right)^2 + \left(\frac{\partial\psi}{\partial x}\right)^2 \right)\, d\Omega.
\end{equation}
In addition, to evaluate the ability of the LES model to reproduce the distribution of energy across spatial scales, we compare the KE spectrum at $t=20$ for each model. %\lander{I am using 20 as a placeholder. This can be decided later when we have the DNS data as well.}
%\anna{Yes, that's fine.}
To get the spectrum, the Fourier transform of each velocity component at the selected time is computed.
Physical wavenumbers $\boldsymbol{k}$ are defined from the grid spacing as $k_{i} = \frac{2\pi}{L}$%\sachin{should it not be $k_{i} = \frac{2\pi}{L}$, where $L$ is the difference of minimum and maximum element of spatial coordinate vector}
, with $i=1, 2$ denoting the
velocity component, so that each mode corresponds to a true spatial frequency. 
The spectral energy density can be obtained as $E(\boldsymbol{k}) = \frac{1}{2}(\sum_{i=1}^2|\hat{u}_{i}|^2)$, where $\hat{u}_i$ indicates the Fourier transform of the $i$-th velocity component.
An isotropic energy spectrum is finally constructed by averaging $E(\boldsymbol{k})$ over annuli of constant wavenumber magnitude $k = |\boldsymbol{k}|$.

Lastly, in Sec.~\ref{sec:mediterranean}, we explore the performance of the different LES models in a realistic geometry, i.e., the Mediterranean sea.

\subsection{The linear and nonlinear BV-Bardina models} \label{sec:bv-bardina}

In this section, we present the numerical results given by the BV-Bardina model with indicator functions $a_{lin}, a_{Leith},$ and $a_S$.

\begin{figure}[htb!]
    \centering
    \begin{tabular}{cccccc}
        \hspace{-0.4cm}& \hspace{-0.4cm}DNS & \hspace{-0.4cm}QGE (coarse) & \hspace{-0.4cm}$a_{lin}$ & \hspace{-0.4cm}$a_{Leith}$ & \hspace{-0.4cm}$a_S$ \\
        \hspace{-0.4cm}$\tilde{\psi}$  & \hspace{-0.4cm}\includegraphics[align=c,scale = 0.25]{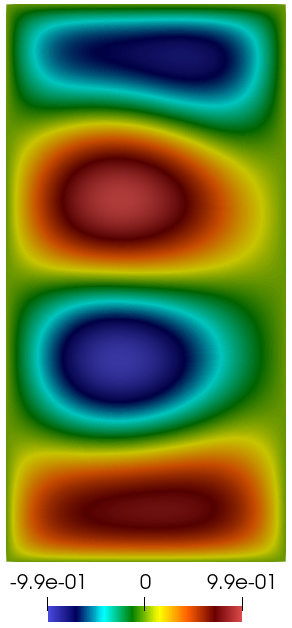} & \hspace{-0.4cm}\includegraphics[align=c,scale = 0.25]{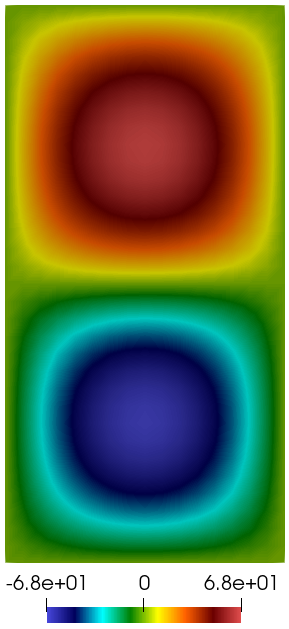} & \hspace{-0.4cm}\includegraphics[align=c,scale = 0.25]{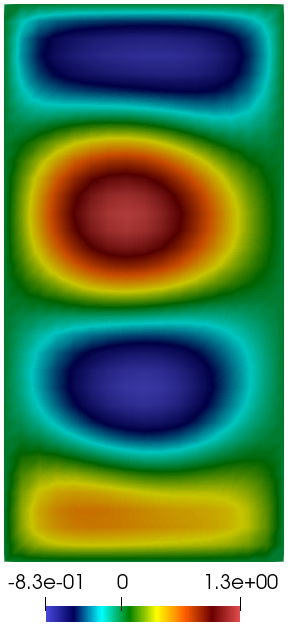} & \hspace{-0.4cm}\includegraphics[align=c,scale = 0.25]{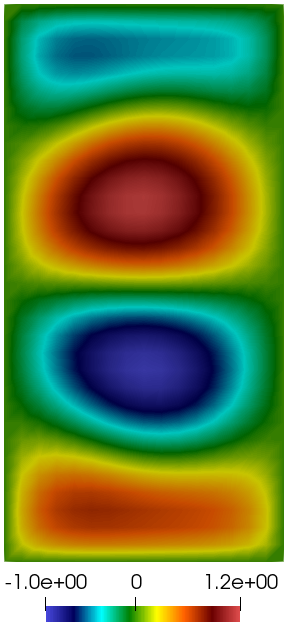} & \hspace{-0.4cm}\includegraphics[align=c,scale = 0.25]{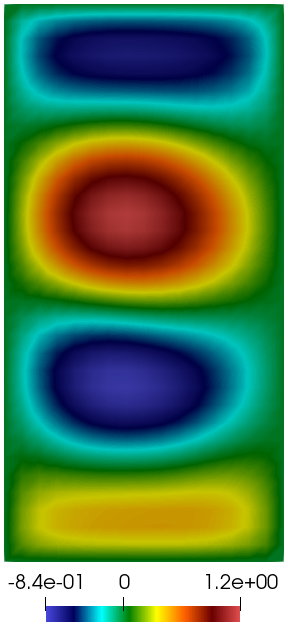} \\
        \hspace{-0.4cm}$\tilde{q}$  & \hspace{-0.4cm}\includegraphics[align=c,scale = 0.25]{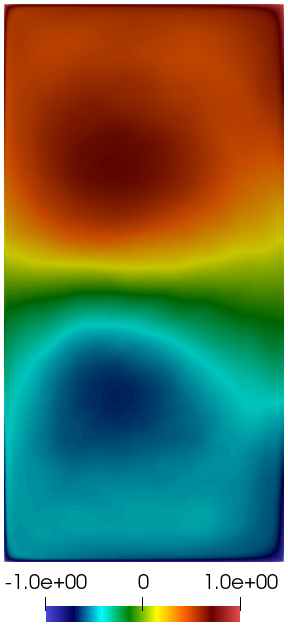} & \hspace{-0.4cm}\includegraphics[align=c,scale = 0.25]{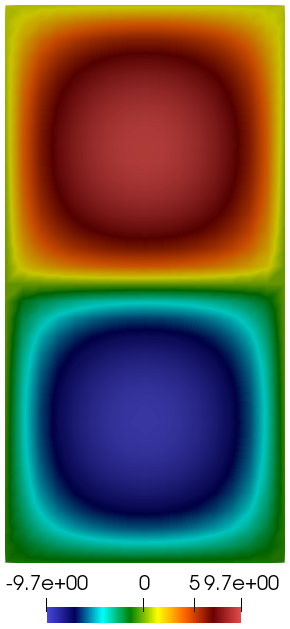} & \hspace{-0.4cm}\includegraphics[align=c,scale = 0.25]{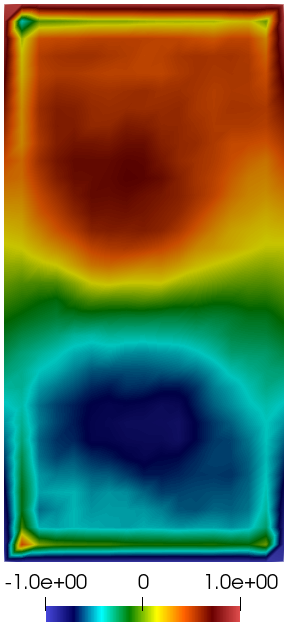} & \hspace{-0.4cm}\includegraphics[align=c,scale = 0.25]{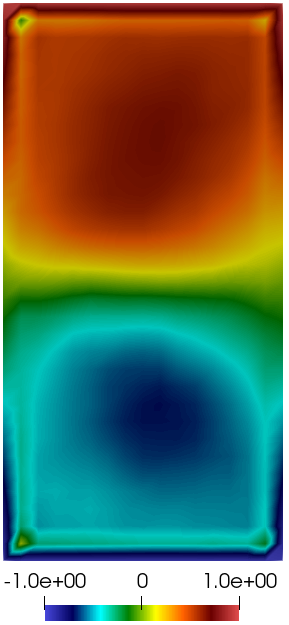} & \hspace{-0.4cm}\includegraphics[align=c,scale = 0.25]{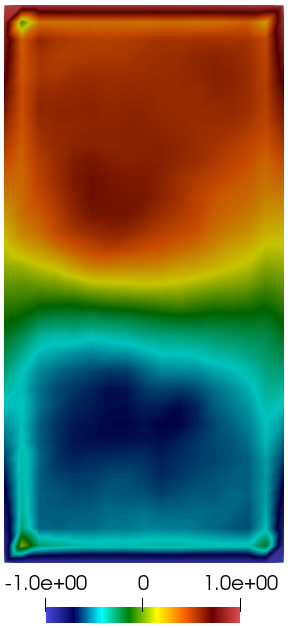}
    \end{tabular}
    \caption{BV-Bardina, experiment 1: Time-averaged stream function $\tilde{\psi}$ (first row) and potential vorticity $\tilde{q}$ (second row) computed from the DNS (first column), QGE with no SGS model (second column), BV-Bardina model with indicator functions $a_{lin}$ (third column), $a_{Leith}$ (fourth column), and $a_S$ (fifth column). The coarse mesh size is $h = 1/16$ and we set $\alpha = h$.}
    \label{fig:Mean_16}
\end{figure}

\begin{figure}[htb!]
    \centering
    \begin{tabular}{ccccc}
        \hspace{-0.4cm}& \hspace{-0.4cm}QGE (coarse) & \hspace{-0.4cm}$a_{lin}$ & \hspace{-0.4cm}$a_{Leith}$ & \hspace{-0.4cm}$a_S$ \\
        \hspace{-0.4cm}$\tilde{q}$ & \hspace{-0.4cm}\includegraphics[align=c,scale = 0.25]{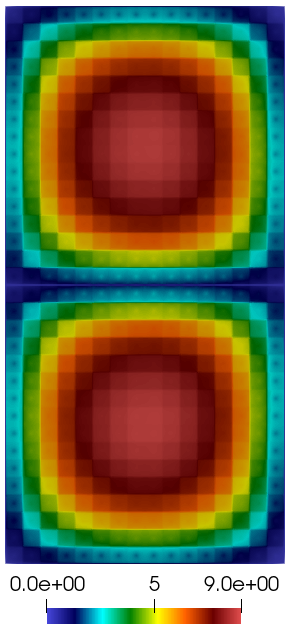} & \hspace{-0.4cm}\includegraphics[align=c,scale = 0.25]{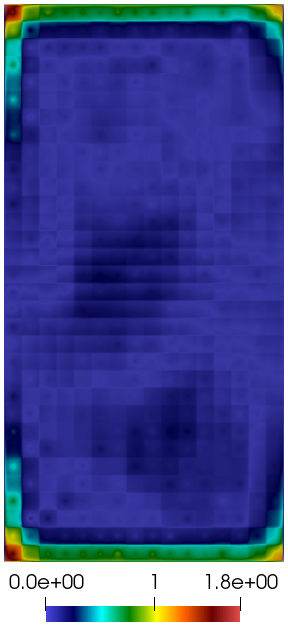} & \hspace{-0.4cm}\includegraphics[align=c,scale = 0.25]{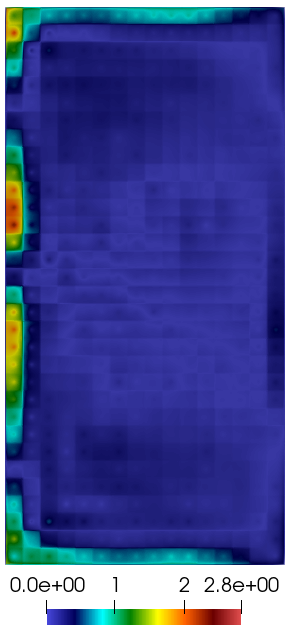} & \hspace{-0.4cm}\includegraphics[align=c,scale = 0.25]{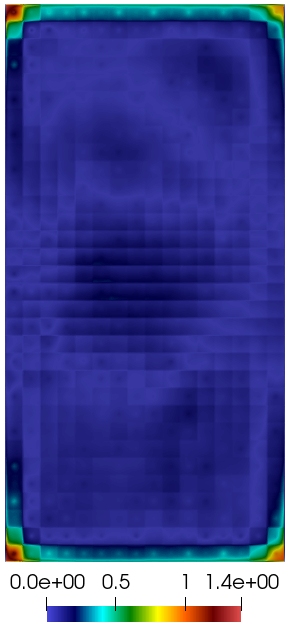}
    \end{tabular}
    \caption{BV-Bardina, experiment 1: Absolute error between $\tilde{q}$ given by the DNS and $\tilde{q}$ computed by QGE with no SGS (first column) and the BV-Bardina model with $a_{lin}$ (second column), $a_{Leith}$ (third column), and $a_S$ (fourth column).
    The coarse mesh size is $h = 1/16$ and we set $\alpha = h$.}
    \label{fig:psiMean_diff_case1}
\end{figure}

We first look at the results given by the BV-Bardina model for experiment 1 with mesh size $h = 1/16$ and $\alpha = h$.
Fig.~\ref{fig:Mean_16}
compares the time averaged stream function $\tilde{\psi}$ and potential vorticity $\tilde{q}$ computed by the DNS (the first column) and the QGE without and with filtering.
As expected \cite{Girfoglio_JCAM2023}, 
$\tilde{\psi}$ given by 
the DNS shows four gyres, while 
$\tilde{\psi}$ computed by the QGE with no SGS model fails to capture the four-gyre pattern and the magnitude of its peaks is about 70 times larger than that of the DNS. 
Compare the first and second panels on the first row of Fig.~\ref{fig:Mean_16}. 
From Fig.~\ref{fig:Mean_16}, we see that BV-Bardina model 
with any of the above mentioned indicator functions recovers the four-gyre pattern
and a magnitude range comparable with the DNS.
The Smagorinsky-like indicator function $a_S$ appears to provide the closest qualitative agreement with the DNS in terms of $\tilde{\psi}$, but the Leith-like model gives a minimum and maximum values of $\tilde{\psi}$ which are closest to the values given by the DNS.
Similarly, 
the transition layer (green region) 
and magnitude range in the time-averaged potential voriticity $\tilde{q}$
are more accurately captured with any BV-Bardina model.
However, this kind of model struggles to recover $\tilde{q}$ in a narrow region close to the boundary, 
as it can be further seen in 
Fig.~\ref{fig:psiMean_diff_case1}, which shows
the absolute errors.

%From hereon, we will refer to the DNS-computed solution as the ``true" solution.

\begin{figure}[htb!]
    \centering
    \begin{tabular}{cccc}
        \hspace{-0.4cm}DNS &\hspace{-0.4cm}$h = 1/16$ & \hspace{-0.4cm}$h = 1/32$ & \hspace{-0.4cm}$h = 1/64$ \\
        \hspace{-0.4cm}\includegraphics[align=c,scale = 0.25]{Exp1_figs/QGE/qMean_DNS.png} & \hspace{-0.4cm}\includegraphics[align=c,scale = 0.25]{Exp1_figs/BVB/16_h/smag/qMean_smag16.png} & \hspace{-0.4cm}\includegraphics[align=c,scale = 0.25]{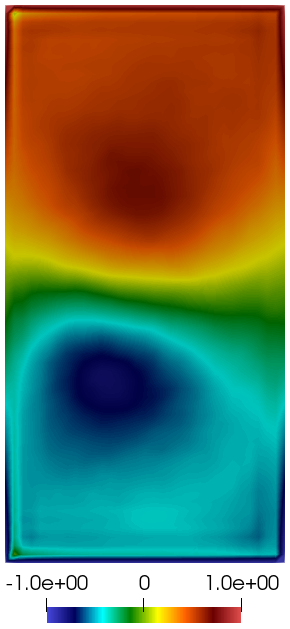} & \hspace{-0.4cm}\includegraphics[align=c,scale = 0.25]{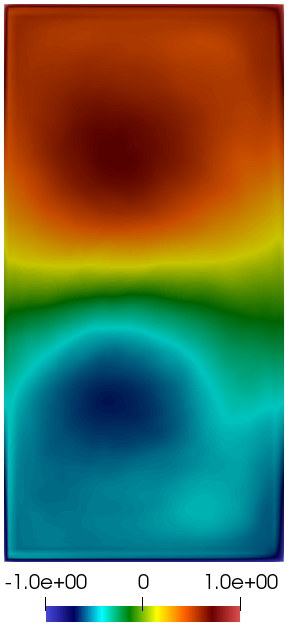}
    \end{tabular}
    \caption{BV-Bardina, experiment 1: Time-averaged potential vorticity $\tilde{q}$ computed by the DNS (leftmost panel) and the BV-Bardina model with indicator function $a_S$ for different mesh sizes $h$. We set $\alpha = h$.}
    \label{fig:qMean_bvb_smag_mesh_refine}
\end{figure}

This instability near the boundary is 
due to a sharp transition required to satisfy boundary condition \eqref{eq:BV5_comp2}. %On very coarse meshes, this is shown as a significant jump near the boundary, leading to relatively large local errors. Consequently, BV-Bardina is less effective in this region on coarse meshes.
As the mesh gets refined, 
the instability disappears.  See Fig.~\ref{fig:qMean_bvb_smag_mesh_refine}, which compares $\tilde{q}$ computed by the BV-Bardina model with $a_S$ for mesh sizes $h = 1/16, 1/32, 1/64$, using $\alpha = h$.

\begin{figure}[htb!]
    \centering
    \begin{tabular}{cc}
        $\tilde{a}_{Leith}, h = 1/16$ & $\tilde{a}_{Leith}, h = 1/32$ \\
        %& %$\tilde{a}_S, h = 1/16$ & $\tilde{a}_S, h = 1/32$ \\
        \includegraphics[align=c,scale = 0.3]{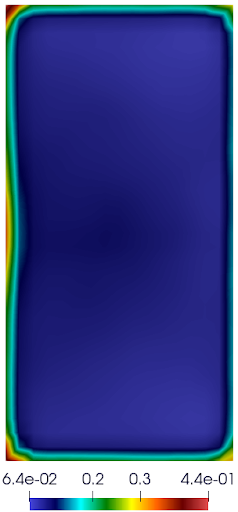} & \includegraphics[align=c,scale = 0.3]{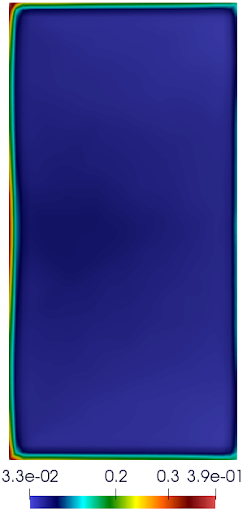} 
        %& 
        % \includegraphics[align=c,scale = 0.3]{Exp1_figs/BVB/a_S.png} & \includegraphics[align=c,scale = 0.3]{Exp1_figs/BVB/a_S32.png}
    \end{tabular}
    \quad \quad
    \begin{tabular}{cc}
        %$\tilde{a}_{Leith}, h = 1/16$ & $\tilde{a}_{Leith}, h = 1/32$ & 
        $\tilde{a}_S, h = 1/16$ & $\tilde{a}_S, h = 1/32$ \\
        \includegraphics[align=c,scale = 0.3]{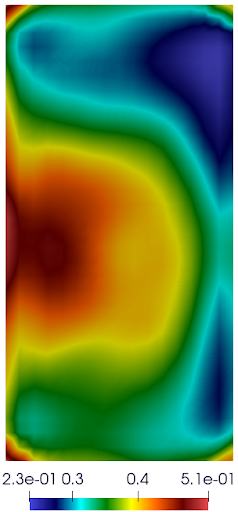} & \includegraphics[align=c,scale = 0.3]{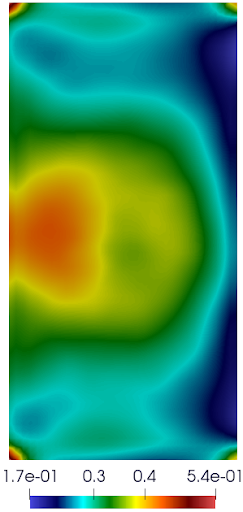}
    \end{tabular}
    \caption{BV-Bardina, experiment 1: Time-averaged indicator function of the Leith-like model (left 2 panels) and the Smagorinsky-like model (right 2 panels) for mesh sizes $h = 1/16, 1/32$. The filtering radius is set to $\alpha = h$.}
    \label{fig:indicator}
\end{figure}

Before presenting the time-dependent quantities, we take a look at the indicator functions. 
Fig.~\ref{fig:indicator} shows the time-averaged indicator function for the Leith-like model and the Smagorinsky-like model for mesh sizes $h = 1/16,1/32$ and filtering radius $\alpha = h$. Both indicator functions take large values in the Northern, Western, and Southern boundaries. However, while the interior values of $\tilde{a}_{Leith}$ are one order of magnitude smaller than the values near the boundary, $\tilde{a}_S$ takes fairly large values in the interior as well, especially towards the center of the basin and near the Western boundary. We note that despite the high selectiveness of the Leith-like indicator, it is sufficient to recover the four gyre structure even a vary coarse mesh like mesh $h = 1/16$, as seen in Fig.~\ref{fig:Mean_16}.

%We can further validate this with the absolute error fields shown in Fig.~\ref{fig:psiMean_diff_case1}. While the coarse QGE solution exhibits relatively large and spatially widespread errors, the BV-Bardina models significantly reduce the magnitude of the absolute difference for $\tilde{q}$. Note that the errors are primarily contained around the boundary. \lander{Consequently, one practical approach is to perform simulate with BV-Bardina on a computational domain larger than that region of interest to minimize the error on the solution around the target region. Lander - not sure if we should add this comment. Need Anna's feedback :)}

\begin{figure}[htb!]
    \centering
    \begin{subfigure}{0.45\linewidth}
        \includegraphics[width = \linewidth]{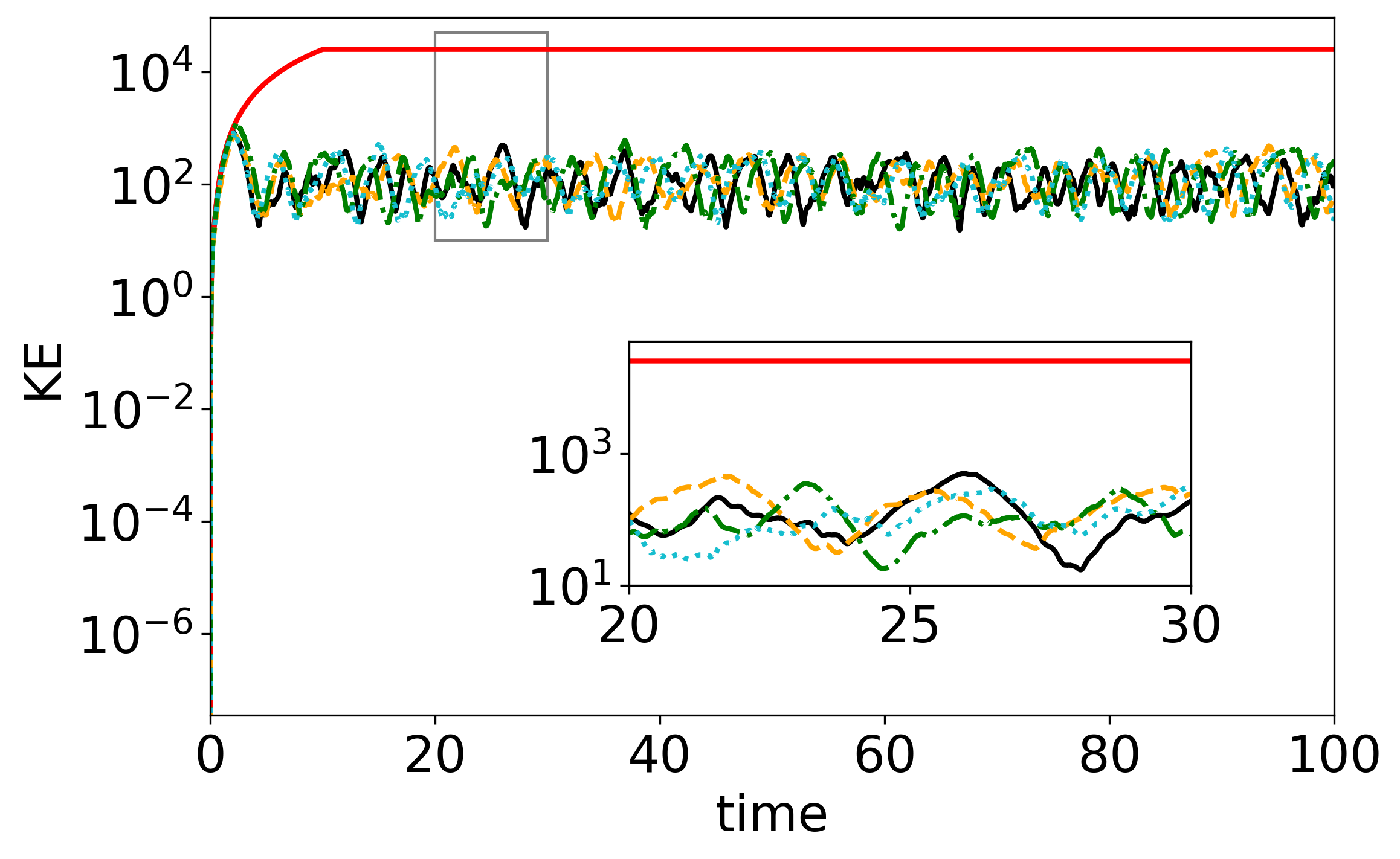}
        \caption{$h = 1/16$}
    \end{subfigure}
    \begin{subfigure}{0.45\linewidth}
        \includegraphics[width = \linewidth]{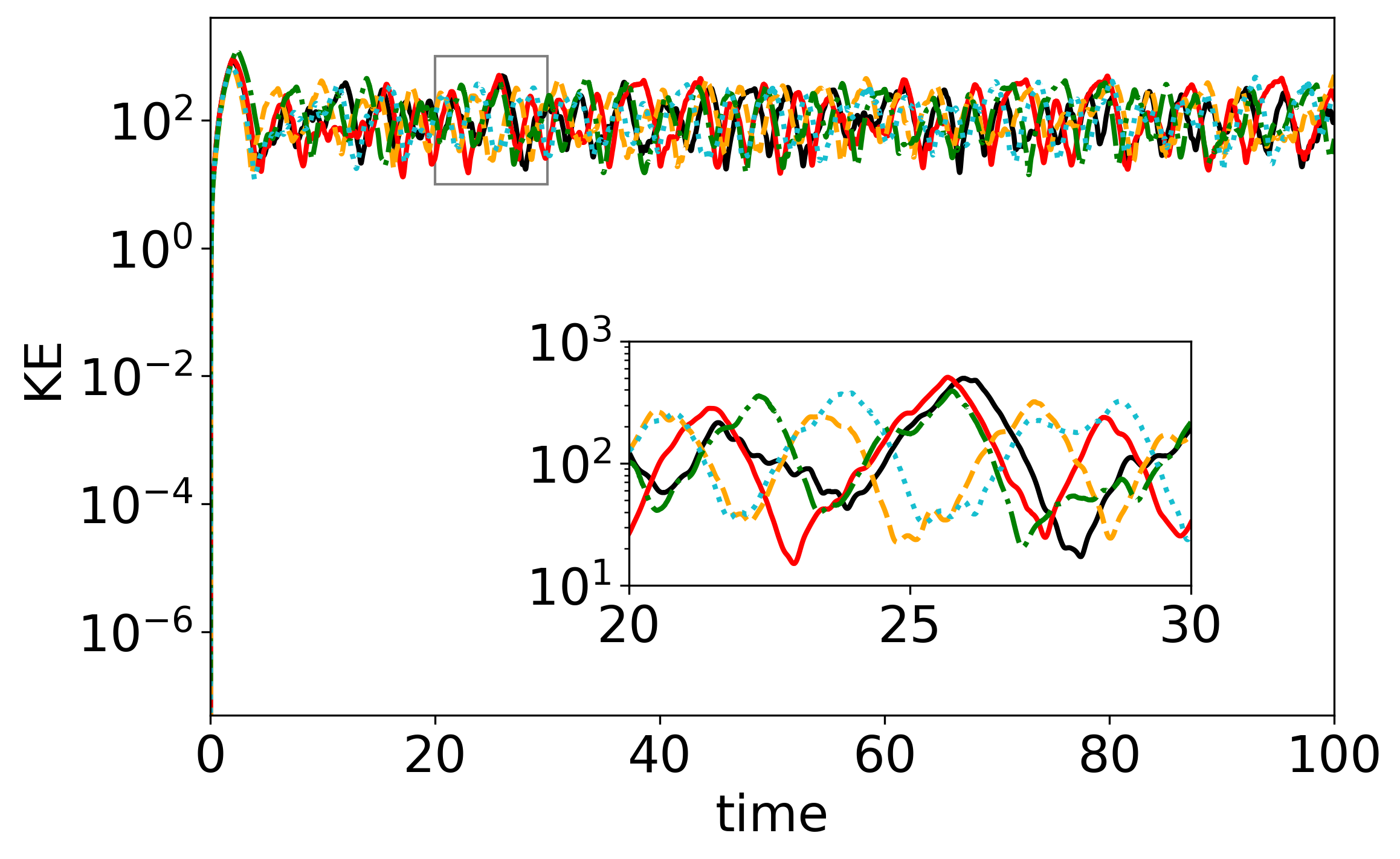}
        \caption{$h = 1/32$}
    \end{subfigure}
    \begin{subfigure}{0.7\linewidth}
        \includegraphics[width = \linewidth]{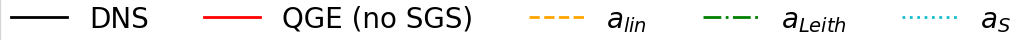}
    \end{subfigure}
    \caption{BV-Bardina, experiment 1: KE in time for the DNS (with mesh $h = 1/256$), QGE with no SGS model,  and BV-Bardina model with indicator functions $a_{lin}$, $a_{Leith}$, and $a_S$ for mesh (a) $h = 1/16$ and (b) $h = 1/32$. We set $\alpha = h$.}
    \label{fig:ke_evol}
\end{figure}

The time evolution of the KE computed by the different models is shown in Fig.~\ref{fig:ke_evol} for meshes $h = 1/16, 1/32$. As expected, with a very coarse mesh (i.e., $h = 1/16$) and no regularization, the KE departs significantly from the DNS reference. Compare the 
red curve (QGE with no SGS) and the 
black curve (DNS) in Fig.~\ref{fig:ke_evol} (a). 
In contrast, the BV-Bardina model with any indicator function introduces sufficient subgrid-scale dissipation to maintain a bounded and physically realistic energy level throughout the simulation. The three indicator functions produce very similar KE evolutions, suggesting that the overall performance of the model is relatively insensitive to the specific choice of the indicator in this regime. Although the BV-Bardina solutions do not remain perfectly in phase with the DNS over the entire time interval, they accurately reproduce both the average kinetic energy level and the amplitude of the low-frequency oscillations.
A comparison between Figs.~\ref{fig:ke_evol} (a)
and \ref{fig:ke_evol} (b) highlights the robustness of the BV-Bardina model with respect to mesh refinement, making it a reliable closure strategy for under-resolved QGE simulations across different resolutions.

With regard to Fig.~\ref{fig:ke_evol}, 
we point out that the BV-Bardina model with the Leith-like indicator is not stable for mesh size $h = 1/32$ and filtering radius $\alpha = h$ unless subiterations are performed as per Remark~\ref{rem:subit}.
This is possibly due to the above-mentioned selectiveness of $a_{Leith}$. Luckily, with even with a rather small stopping tolerance in \eqref{eq:stop} like 
$\epsilon = 10^{-6}$, only an average of three iterations per time step are needed.  
In contrast, 
the BV-Bardina model with the Smagorinsky-like indicator
never requires more than one iteration per time step, indicating robustness with respect to mesh size changes, 
as already observed in \cite{Maulik2016}.
%\lander{, which is in contrast to the conclusions made for the BV-$\alpha$ models (Lander - again, not sure if we will include this)}.

\begin{figure}[htb!]
    \centering
    \begin{subfigure}{0.45\linewidth}
        \includegraphics[width = \linewidth]{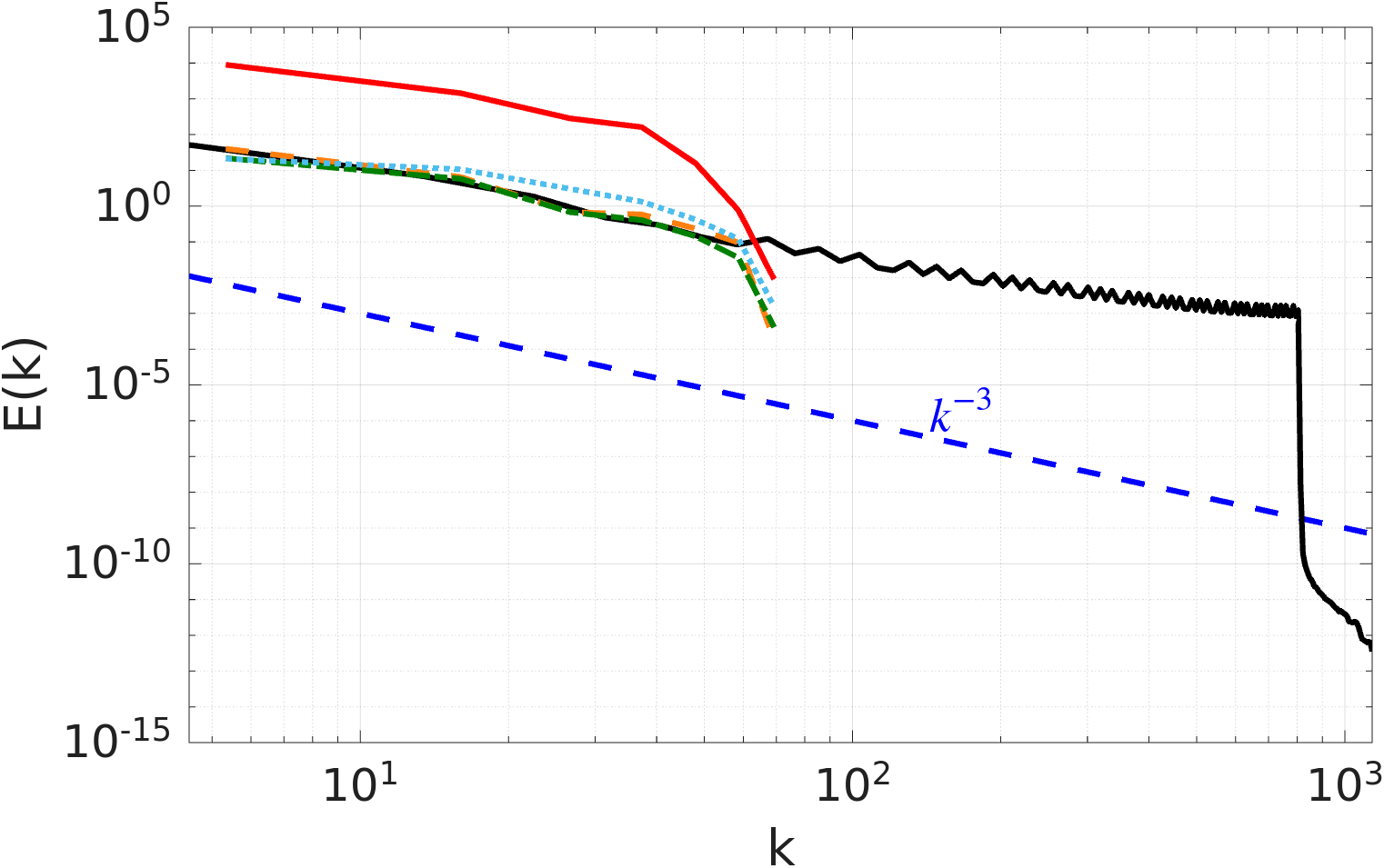}
        \caption{$h = 1/16$}
    \end{subfigure}
    \begin{subfigure}{0.45\linewidth}
        \includegraphics[width = \linewidth]{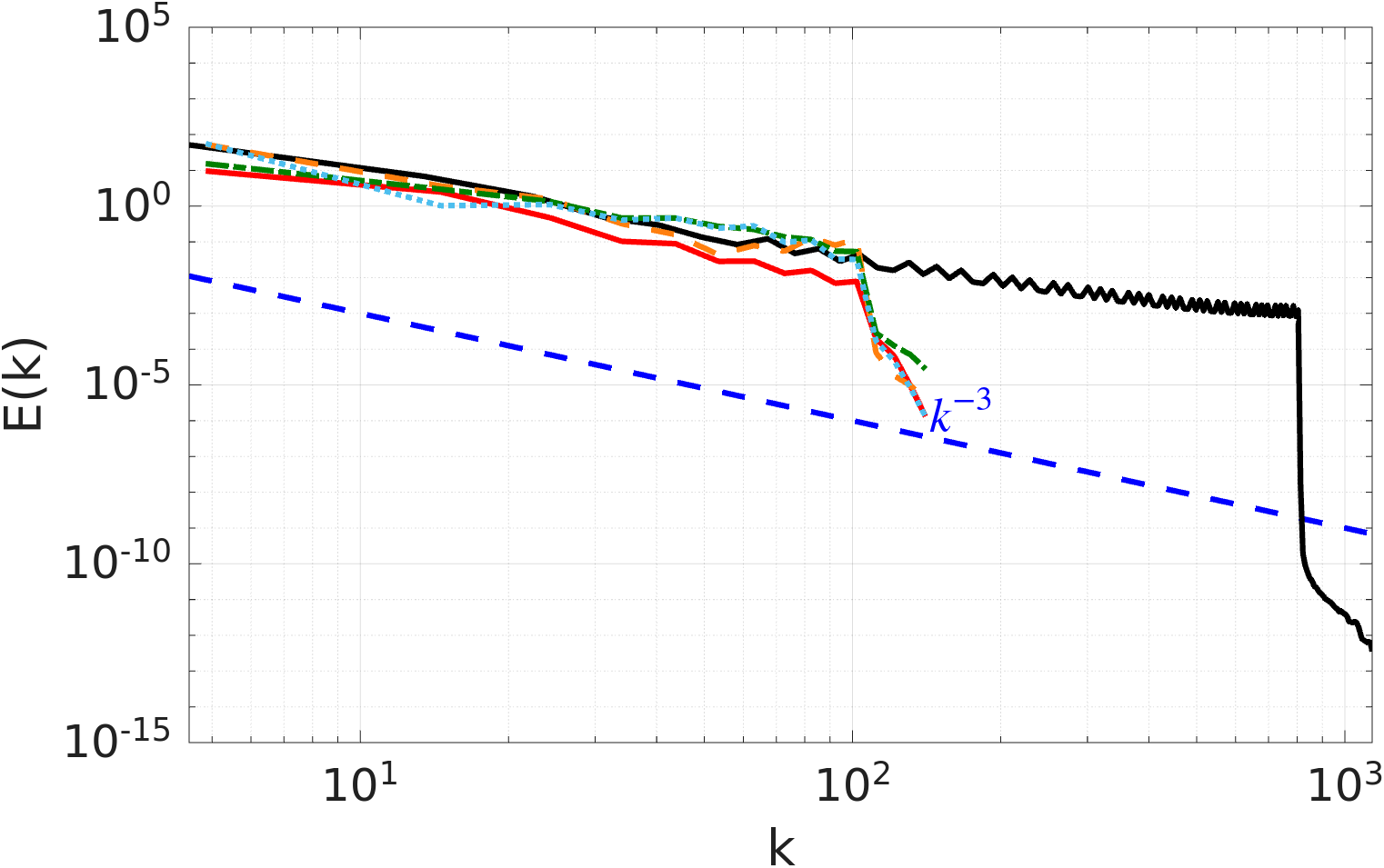}
        \caption{$h = 1/32$}
    \end{subfigure}
    \begin{subfigure}{0.7\linewidth}
        \includegraphics[width = \linewidth]{Exp1_figs/BVB/legend.png}
    \end{subfigure}
    \caption{BV-Bardina, experiment 1: Spectrum of the KE given by the DNS 
    compared with the spectrum given by the
    QGE with no SGS model,  and the BV-Bardina model with indicator functions $a_{lin}$, $a_{Leith}$, and $a_S$ at time $t = 20$ for mesh (a) $h = 1/16$ and (b) $h = 1/32$. We set $\alpha = h$. %\anna{Minor comment: the green line seems solid to me, no?}\sachin{I fixed it.}
    %\anna{Can we make the colored lines thicker like in Fig. 6 so that they're easier to see?}
    %\sachin{updated to size 3, earlier it was 2, in Fig. 6 we are using 2 but the Kinetic energy figures are generated by python and spectrum figures are generated by Matlab that's why they might look a bit different}
    }
    \label{fig:ke_spec_16}
\end{figure}

Fig.~\ref{fig:ke_spec_16} compares the KE spectra for all the models under consideration at $t=20$ using coarse meshes $h = 1/16, 1/32$. The QGE simulation with mesh $h = 1/16$ exhibits a substantial excess of energy across the resolved wavenumber range, indicating an insufficient dissipation due to lack of model for the unresolved scales. In contrast, the BV-Bardina model with any indicator function significantly reduces the energy content and produces spectra that are much closer to the DNS spectrum, all exhibiting the correct slope over the resolved scales. 
The results for mesh  
$h = 1/32$ confirm that the BV-Bardina spectra reproduce the DNS trend well at the large and intermediate scales, while the largest discrepancies remain confined to the highest resolved wavenumbers, where the effects of under-resolution are expected to be strongest. The differences among the three indicator functions are relatively small compared to the discrepancy between the coarse QGE and DNS solutions, suggesting that the introduction of the filtering mechanism itself plays the dominant role in recovering the correct energy distribution. 
The results also indicate that mesh refinement from $h=1/16$ to $h=1/32$ improves agreement with the DNS, particularly at intermediate and high wavenumbers.

The results of experiment 1 indicate that mesh size $h = 1/16$ is too coarse to accurately capture the time-averaged potential vorticity $\tilde{q}$, particularly near the boundary, while the mesh size $h = 1/32$ provides a substantially more accurate reconstruction of the flow. Hence, from now on we will consider only mesh size $h = 1/32$, which is eight times the mesh size used for the DNS.

\begin{figure}[htb!]
    \centering
    \begin{tabular}{cccccc}
        \hspace{-0.4cm}& \hspace{-0.4cm}DNS & \hspace{-0.4cm}QGE (coarse) & \hspace{-0.4cm}$a_{lin}$ & \hspace{-0.4cm}$a_{Leith}$ & \hspace{-0.4cm}$a_S$ \\
        \hspace{-0.4cm}$\tilde{\psi}$  & \hspace{-0.4cm}\includegraphics[align=c,scale = 0.25]{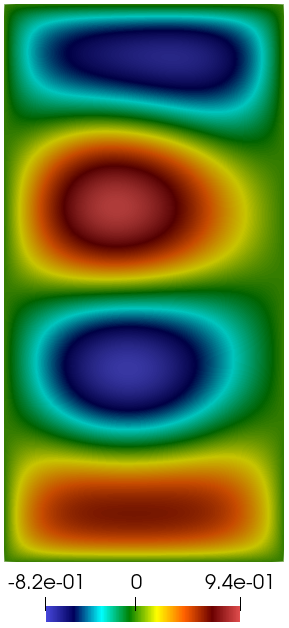} & \hspace{-0.4cm}\includegraphics[align=c,scale = 0.25]{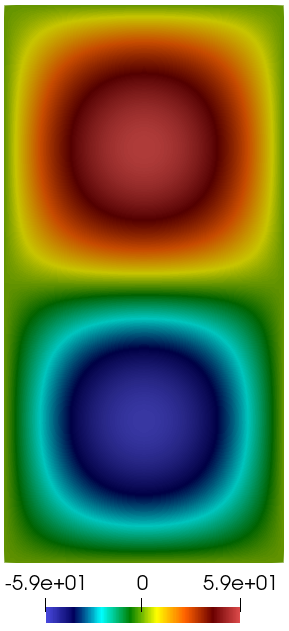} & \hspace{-0.4cm}\includegraphics[align=c,scale = 0.25]{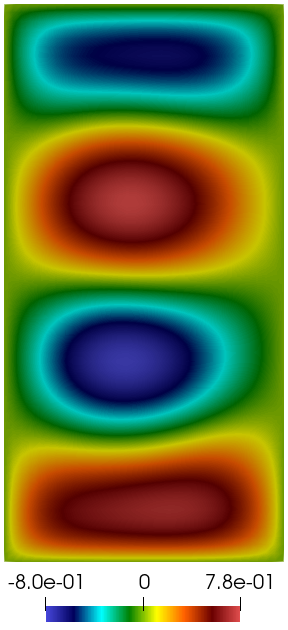} & \hspace{-0.4cm}\includegraphics[align=c,scale = 0.25]{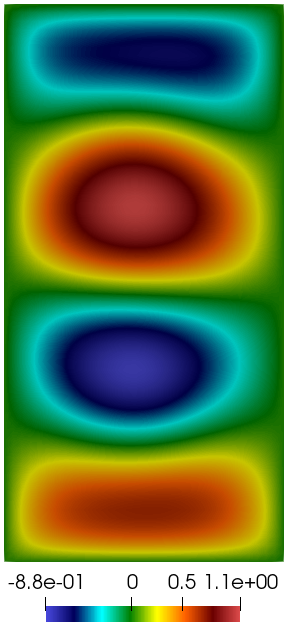} & \hspace{-0.4cm}\includegraphics[align=c,scale = 0.25]{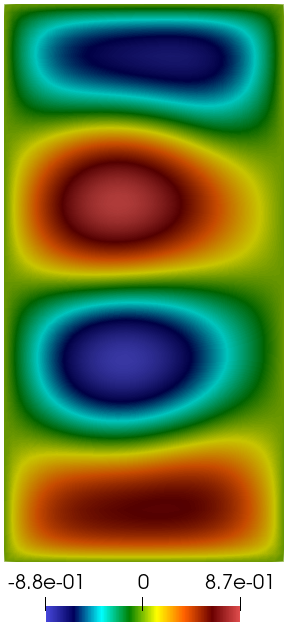} \\
        \hspace{-0.4cm}$\tilde{q}$  & \hspace{-0.4cm}\includegraphics[align=c,scale = 0.25]{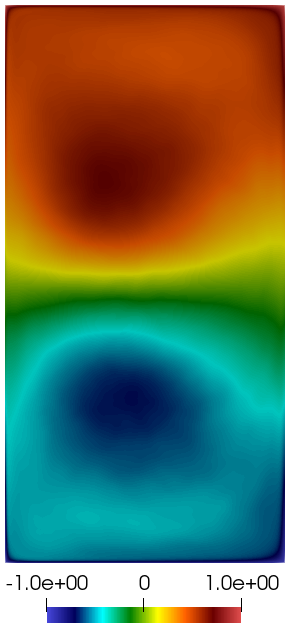} & \hspace{-0.4cm}\includegraphics[align=c,scale = 0.25]{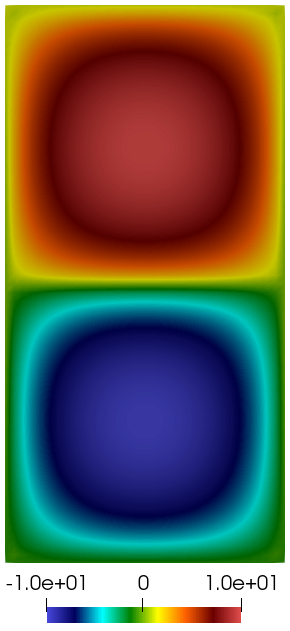} & \hspace{-0.4cm}\includegraphics[align=c,scale = 0.25]{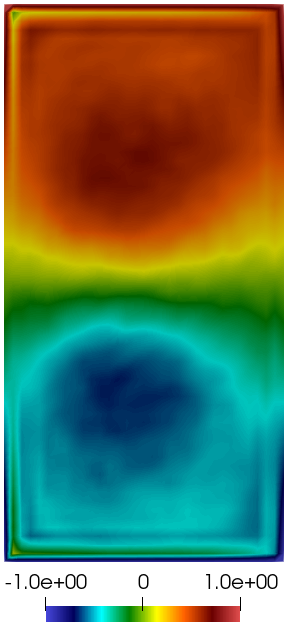} & \hspace{-0.4cm}\includegraphics[align=c,scale = 0.25]{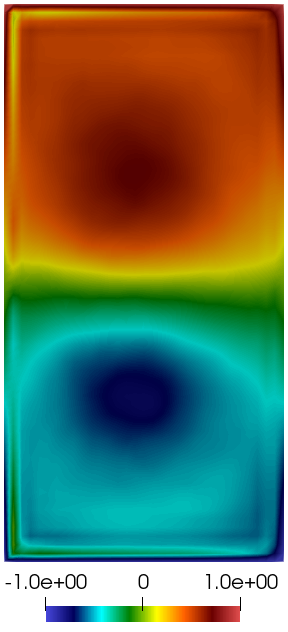} & \hspace{-0.4cm}\includegraphics[align=c,scale = 0.25]{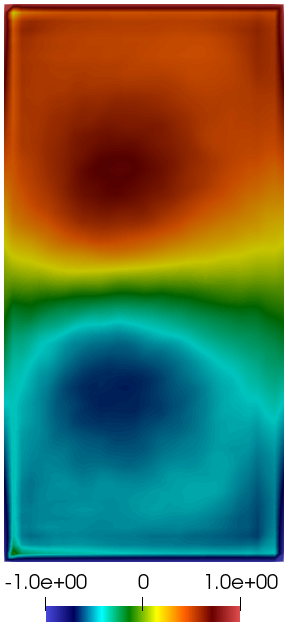}
    \end{tabular}
    \caption{BV-Bardina, experiment 2: Time-averaged stream function $\tilde{\psi}$ (first row) and potential vorticity $\tilde{q}$ (second row) computed from the DNS (first column), QGE with no SGS model (second column), BV-Bardina model with indicator functions $a_{lin}$ (third column), $a_{Leith}$ (fourth column), and $a_S$ (fifth column). The coarse mesh size is $h = 1/32$ and we set $\alpha = h$.}
    \label{fig:Mean_32_exp2}
\end{figure}

Let us turn to experiment 2. 
Fig.~\ref{fig:Mean_32_exp2}
compares the time averaged stream function $\tilde{\psi}$ and potential vorticity $\tilde{q}$ computed by the DNS (the first column) and the QGE without and with filtering.
As for experiment 1, 
$\tilde{\psi}$ given by 
the DNS shows four gyres, while 
$\tilde{\psi}$ computed by the QGE and no SGS model fails to capture the four-gyre pattern.
From 
Fig.~\ref{fig:Mean_32_exp2}, we see that BV-Bardina model 
with any of the considered indicator functions allows us to recover the four-gyre pattern
and a magnitude range comparable with the DNS.
We observe that the BV-Bardina model with $a_{lin}$
and $a_{S}$ over-smooths the positive peaks in $\tilde{\psi}$, while the same model with $a_{Leith}$ does not smooth them sufficiently. This is in line with our previous observation that $a_{Leith}$ is a more selective
indicator function. 
We note that also in the case of experiment 2
the BV-Bardina model with $a_{Leith}$ is unstable without subiterations (see Remark \ref{rem:subit}).
Just like
Figs.~\ref{fig:Mean_16} and \ref{fig:qMean_bvb_smag_mesh_refine}, Fig.~\ref{fig:Mean_32_exp2} shows that the BV-Bardina model struggles to recover $\tilde{q}$ in a narrow region close to the boundary.

%In this setting, coarse QGE with no SGS model and BV-Bardina with $a_{Leith}$ fail to recover the shape of the gyres or the range of values of $\tilde{\psi}$ and $\tilde{q}$. In contrast, BV-Bardina with $a_{lin}$ and $a_S$ substantially improves the accuracy of the coarse simulations. Both models accurately capture the overall gyre structures in the true $\tilde{\psi}$, and have minimum and maximum values which are much closer to that of the DNS. Furthermore, $a_S$ provides a a more accurate approximation of $\tilde{q}$ around the boundary compared to $a_{lin}$.

\begin{figure}[htb!]
    \centering
    \begin{subfigure}{0.45\linewidth}
        \includegraphics[width = \linewidth]{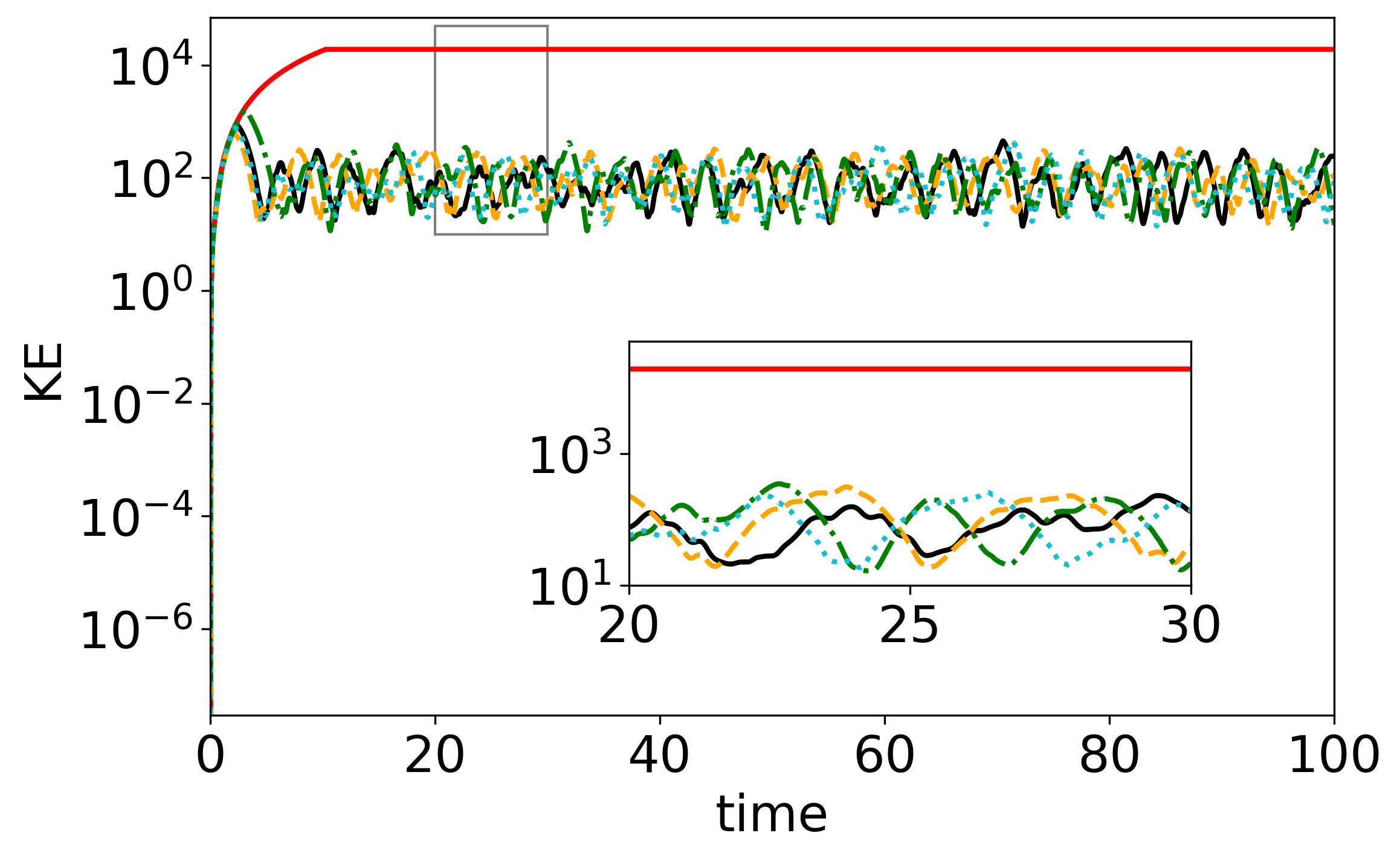} 
        \caption{Kinetic energy}
    \end{subfigure}
    \begin{subfigure}{0.45\linewidth}
        \includegraphics[width = \linewidth]{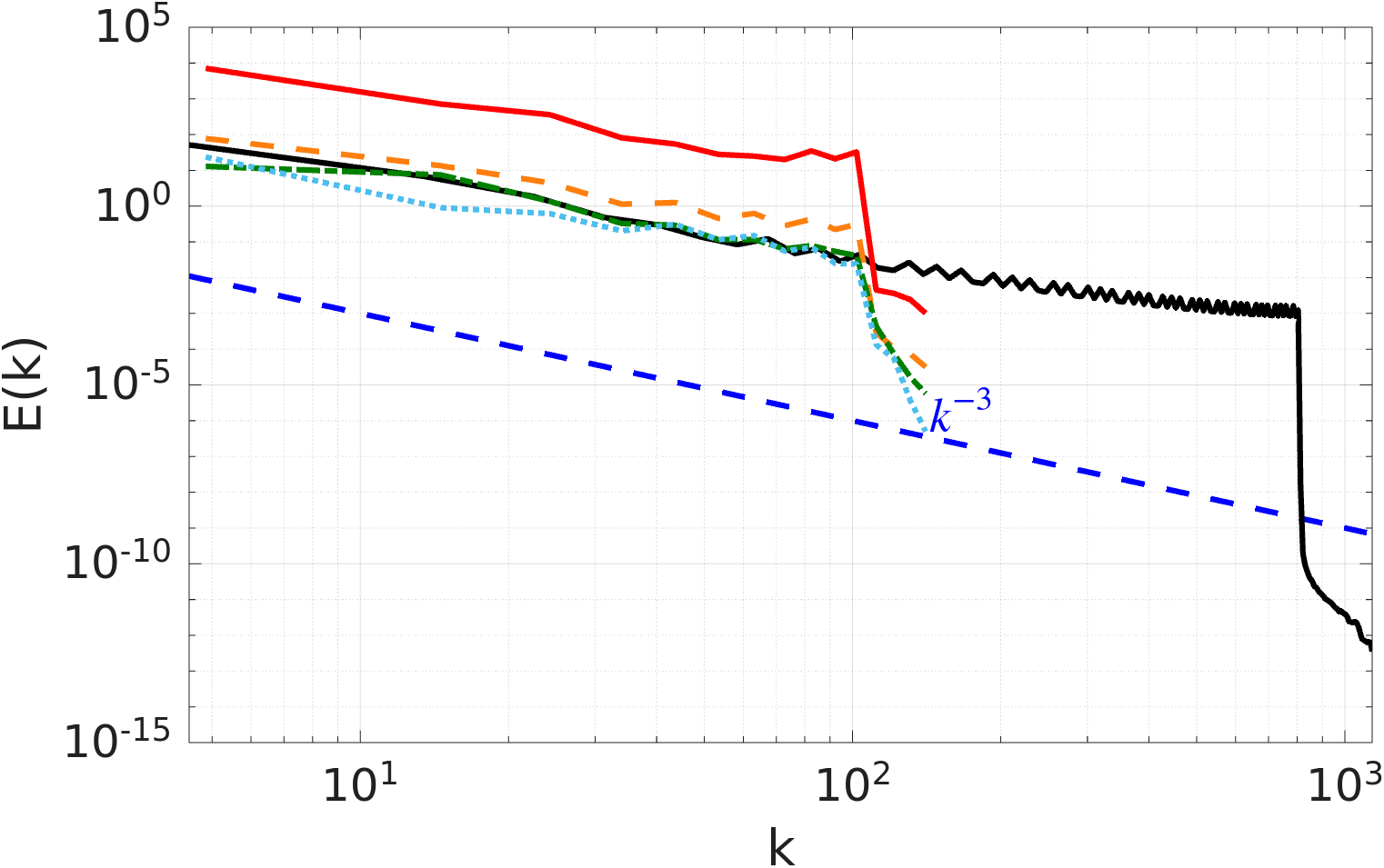}
        \caption{KE spectra}
    \end{subfigure}
    \begin{subfigure}{0.7\linewidth}
        \includegraphics[width = \linewidth]{Exp1_figs/BVB/legend.png}
    \end{subfigure}
    \caption{BV-Bardina, experiment 2: (a) Evolution of the KE and (b) KE spectrum at time $t = 20$ given by the DNS, the
    QGE with no SGS model,  and the BV-Bardina model with indicator functions $a_{lin}$, $a_{Leith}$, and $a_S$ for mesh $h = 1/32$. We set $\alpha = h$. %\anna{On the right panel, the green line is solid, not dash-dotted. And please make the colored lines thicker too, to match the thinkness of the lines on the left.}\sachin{updated.}
    }
    \label{fig:ke_kespec_16_exp2}
\end{figure}

Fig.~\ref{fig:ke_kespec_16_exp2} displays the evolution of the KE and its spectra at $t = 20$
for all the model under consideration, using coarse mesh $h = 1/32$ and $\alpha = h$. As expected from Fig.~\ref{fig:Mean_32_exp2}, the coarse mesh QGE simulation with no SGS model overestimates the KE. On the other hand, the BV-Bardina model with all indicators recovers both the temporal evolution of the kinetic energy and its spectrum over the resolved scales with good accuracy.
Among the tested indicators, $a_{Leith}$  provides the closest agreement with the DNS spectrum, particularly at intermediate and high resolved wavenumbers, suggesting a more accurate representation of the energy transfer across scales at the price of an increased computational cost due to the subiterations. This improved spectral accuracy is consistent with the reduced over-smoothing observed in Fig.~\ref{fig:Mean_32_exp2}.

Overall, in terms of robustness with respect to changes in the mesh size $h$ and physical parameters $Re$ and $Ro$, 
the BV-Bardina model with the Smagorinsky-like indicator function $a_S$ seems to be the most reliable choice. However, for more complex cases, additional dissipation may be required, and the BV-Bardina with $a_{lin}$ can provide a more stable alternative, as will be shown in Sec.~\ref{sec:mediterranean}. %\anna{Add that in more complex cases, one might have to resort to a lin and point to Sec. 4.3.} \lander{I added a sentence.}
%Thus, of the three indicator functions considered in this section we will consider only $a_S$ in the next sections. 

\subsection{The QGE-EFR model} \label{sec:qge-efr}

In this section, we assess the performance of the QGE-EFR model with the indicator functions $a_{lin}$
and $a_{S}$, but replace $a_{Leith}$ with $a_D$.
%Preliminary numerical tests indicate that EFR is unstable on the coarsest mesh $h=1/16$ considered. Moreover, even when stable, its accuracy for $h=1/16$ with $\alpha=h$ and $\alpha=2h$ is unsatisfactory. Therefore, in the following, we present results only for the finer mesh, $h=1/32$. A possible explanation is that, for mesh size $h=1/16$, the amount of artificial dissipation introduced by EFR is insufficient, particularly for $\alpha=h$ and $\alpha=2h$, which results in either numerical instability or poor accuracy.
The reason why we do not consider $a_{Leith}$ any further is because, in the case of the BV-Bardina model, it required subiterations
for stability in both experiment 1 and 2. 

\begin{figure}[htb!]
    \centering
    \begin{tabular}{cccccc}
        \hspace{-0.4cm}& \hspace{-0.4cm}DNS & \hspace{-0.4cm}QGE (coarse) & \hspace{-0.4cm}$a_{lin}$ & \hspace{-0.4cm}$a_{D}$ & \hspace{-0.4cm}$a_S$ \\
        \hspace{-0.4cm}$\tilde{\psi}$  & \hspace{-0.4cm}\includegraphics[align=c,scale = 0.25]{Exp1_figs/QGE/psiMean_DNS.png} & \hspace{-0.4cm}\includegraphics[align=c,scale = 0.25]{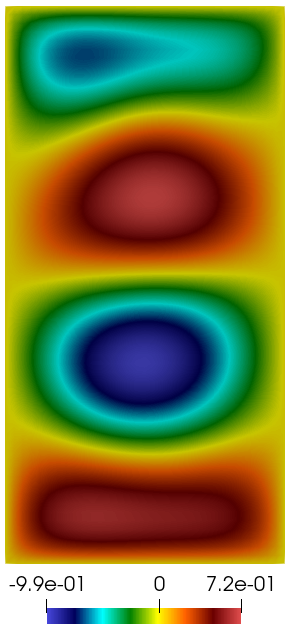} & \hspace{-0.4cm}\includegraphics[align=c,scale = 0.25]{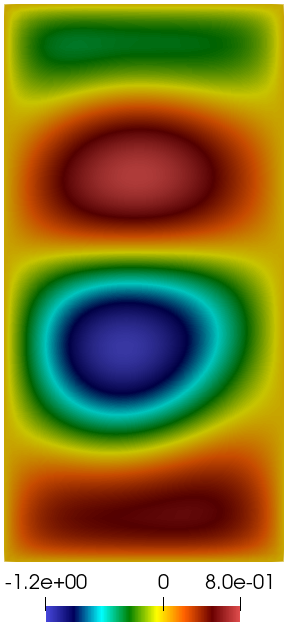} & \hspace{-0.4cm}\includegraphics[align=c,scale = 0.25]{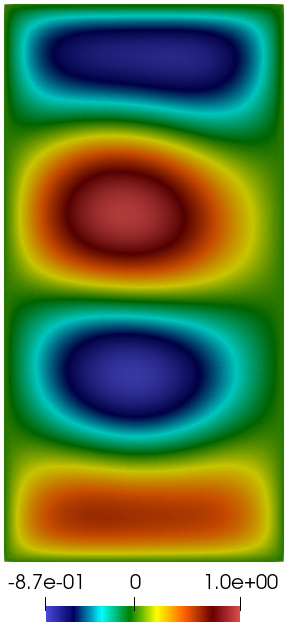} & \hspace{-0.4cm}\includegraphics[align=c,scale = 0.25]{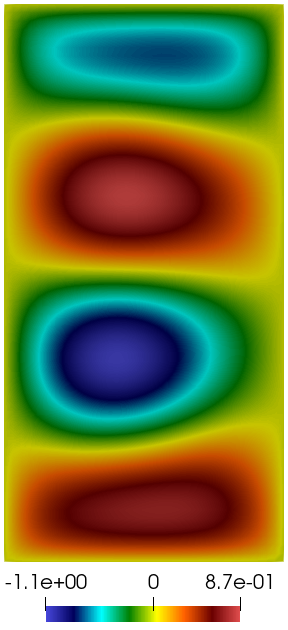} \\
        \hspace{-0.4cm}$\tilde{q}$  & \hspace{-0.4cm}\includegraphics[align=c,scale = 0.25]{Exp1_figs/QGE/qMean_DNS.png} & \hspace{-0.4cm}\includegraphics[align=c,scale = 0.25]{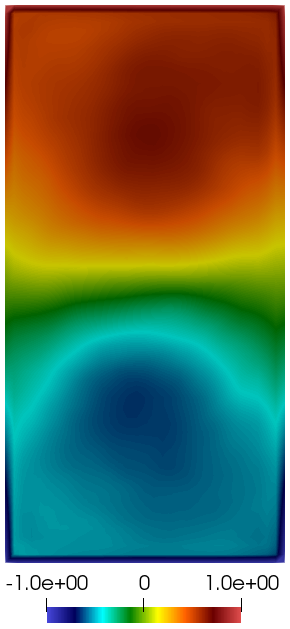} & \hspace{-0.4cm}\includegraphics[align=c,scale = 0.25]{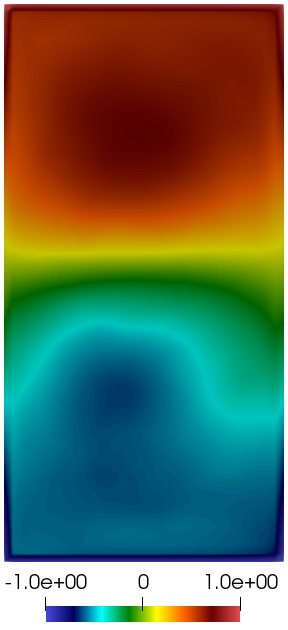} & \hspace{-0.4cm}\includegraphics[align=c,scale = 0.25]{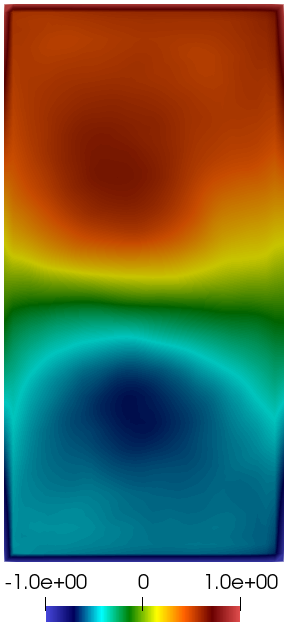} & \hspace{-0.4cm}\includegraphics[align=c,scale = 0.25]{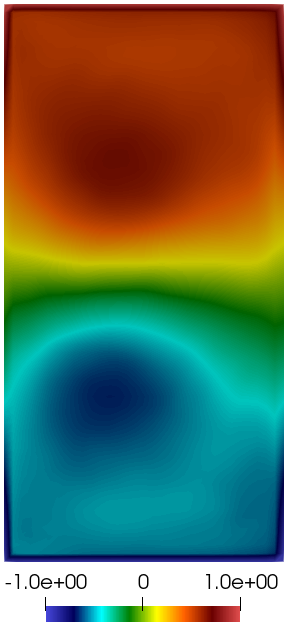}
    \end{tabular}
    \caption{QGE-EFR, experiment 1: Time-averaged stream function $\tilde{\psi}$ (first row) and potential vorticity $\tilde{q}$ (second row) computed from the DNS (first column), QGE with no SGS model (second column), QGE-EFR model with indicator functions $a_{lin}$ (third column), $a_{D}$ (fourth column), and $a_S$ (fifth column). The coarse mesh size is $h = 1/32$ and we set $\alpha = h$.}
    \label{fig:Mean_32_efr}
\end{figure}

Fig.~\ref{fig:Mean_32_efr} presents the time-averaged stream function $\tilde{\psi}$ and potential vorticity $\tilde{q}$ for experiment 1 computed by the DNS, coarse mesh QGE simulation with no filtering, and QGE-EFR with indicator functions $a_{lin}$, $a_D$, and $a_S$.
We consider coarse mesh size $h = 1/32$ and
set to $\alpha = h$. 
All the coarse simulations provide an approximation of $\tilde{\psi}$ with four gyres
and they differ mostly on their reconstruction of the northernmost gyre, where the QGE-EFR model with $a_D$ provides the most accurate shape qualitatively. It also gives the range of values closest to the DNS among all other models. In contrast to the BV-Bardina model, the QGE-EFR model captures well the time-averaged potential vorticity $\tilde{q}$ near the boundaries, as it does not exhibit any spurious oscillations. Compare the bottom row in Fig.~\ref{fig:Mean_16} and Fig.~\ref{fig:qMean_bvb_smag_mesh_refine} with the bottom row in Fig.~\ref{fig:Mean_32_efr}. 

\begin{figure}[htb!]
    \centering
    \begin{subfigure}{0.45\linewidth}
        \includegraphics[width = \linewidth]{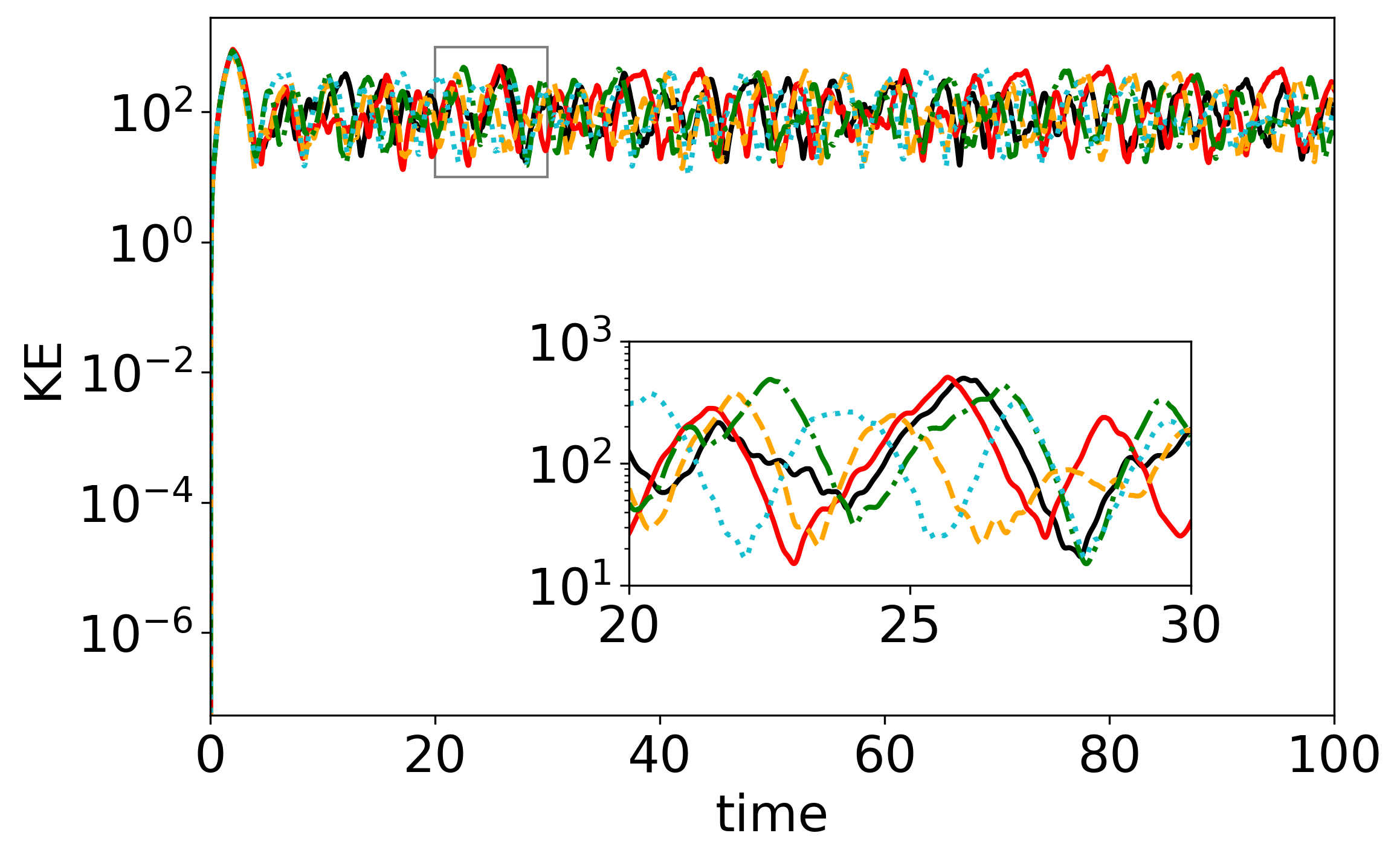}
        \caption{Kinetic energy}
    \end{subfigure}
    \begin{subfigure}{0.45\linewidth}
        \includegraphics[width = \linewidth]{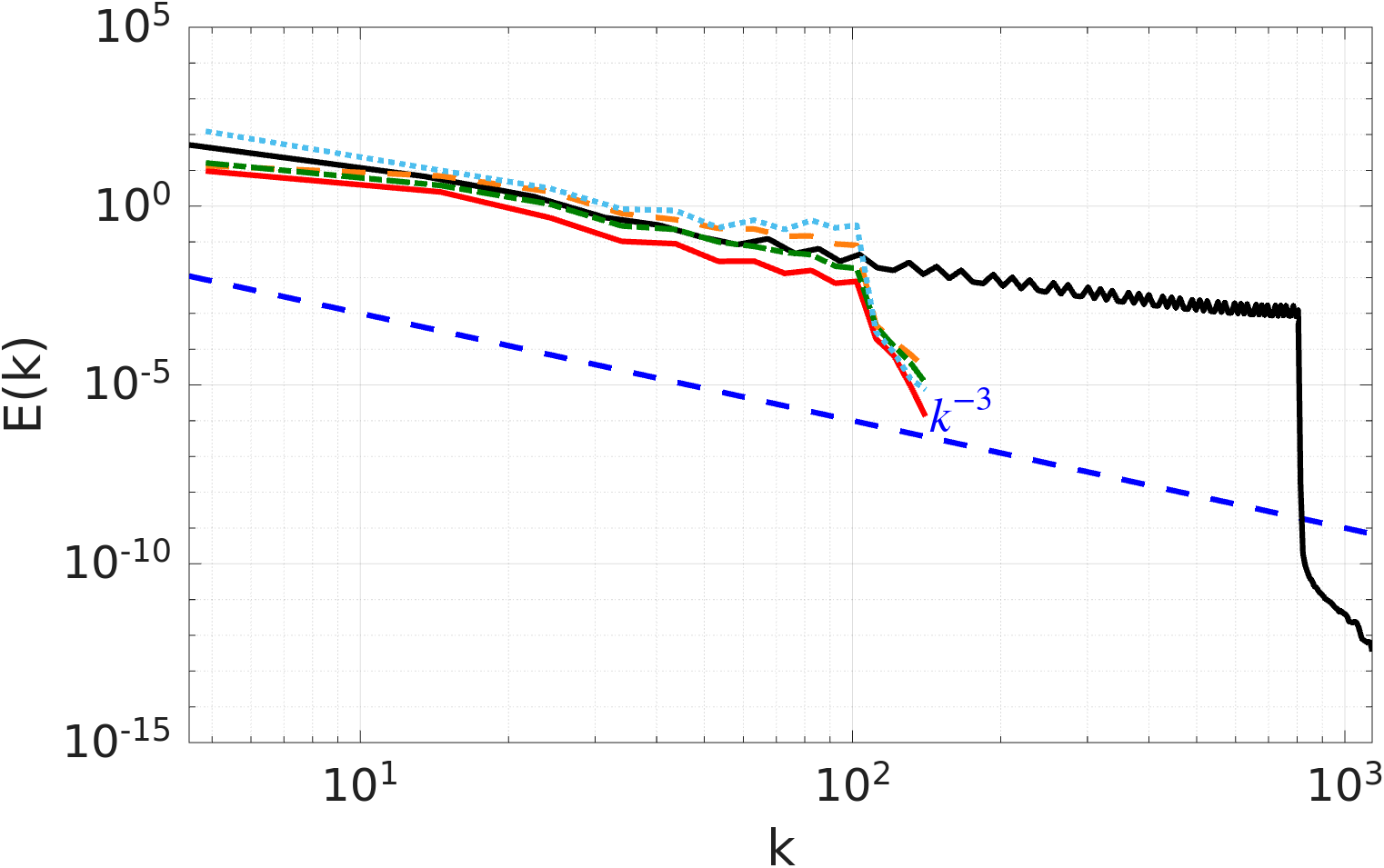}
        \caption{KE spectra}
    \end{subfigure}
    \begin{subfigure}{0.7\linewidth}
        \includegraphics[width = \linewidth]{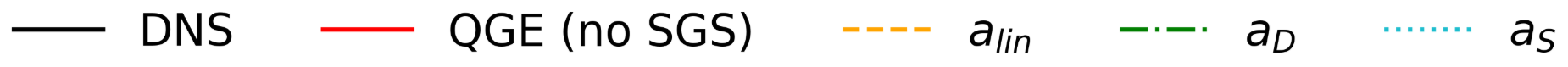}
    \end{subfigure}
    \caption{QGE-EFR, experiment 1:  (a) Evolution of the KE and (b) KE spectrum at time $t = 20$ given by the DNS, the
    QGE with no SGS model,  and the QGE-EFR model with indicator functions $a_{lin}$, $a_{D}$, and $a_S$ for mesh $h = 1/32$. We set $\alpha = h$. %\anna{Can we make the colored lines in the right panel thicker?}\sachin{done.}
    }
    \label{fig:ke_evol_efr}
\end{figure}

Fig.~\ref{fig:ke_evol_efr} shows the KE evolution and spectrum at $t = 20$ for the different QGE-EFR models with coarse mesh $h = 1/32$ and $\alpha = h$. All models that use a coarse mesh produce stable KE evolutions with amplitude comparable to that of the DNS. Although none of the QGE-EFR solutions remains consistently in phase with the DNS, the QGE-EFR model with indicator function $a_{D}$ seems to provide the most accurate reconstruction of the KE. 
Similarly, while all coarse mesh models give a spectrum that compares well with the DNS for the resolved scales, 
the QGE-EFR model with
$a_{D}$ produces the KE spectrum closest to the DNS spectrum. 

% \begin{figure}[htb!]
%     \centering
%     \begin{subfigure}{0.4\linewidth}
%         \includegraphics[width = \linewidth]{Exp1_figs/EFR/16_h/EnergySpectrum_20.png}
%         \caption{$h = 1/16$, $\alpha = h$}
%     \end{subfigure}
%     \begin{subfigure}{0.4\linewidth}
%         \includegraphics[width = \linewidth]{Exp1_figs/EFR/32_h/EnergySpectrum_20.png}
%         \caption{$h = 1/32$, $\alpha = h$}
%     \end{subfigure}
%     \begin{subfigure}{0.7\linewidth}
%         \includegraphics[width = \linewidth]{Exp1_figs/EFR/legend.png}
%     \end{subfigure}
%     \caption{QGE-EFR, experiment 1: Spectrum of the kinetic energy for the DNS (with mesh $h = 1/256$), QGE with no subgrid scale (SGS) model,  and QGE-EFR with indicator functions $a_{lin}$, $a_{D}$, and $a_S$ at time $t = 20$.}
%     \label{fig:ke_spec_efr}
% \end{figure}

\begin{figure}[htb!]
    \centering
    \begin{tabular}{cccccc}
        \hspace{-0.4cm}& \hspace{-0.4cm}DNS & \hspace{-0.4cm}QGE (coarse) & \hspace{-0.4cm}$a_{lin}$ & \hspace{-0.4cm}$a_{D}$ & \hspace{-0.4cm}$a_S$ \\
        \hspace{-0.4cm}$\tilde{\psi}$  & \hspace{-0.4cm}\includegraphics[align=c,scale = 0.25]{Exp2_figs/QGE/psiMean_DNS.png} & \hspace{-0.4cm}\includegraphics[align=c,scale = 0.25]{Exp2_figs/QGE-noSGS/psiMean_qge32.png} & \hspace{-0.4cm}\includegraphics[align=c,scale = 0.25]{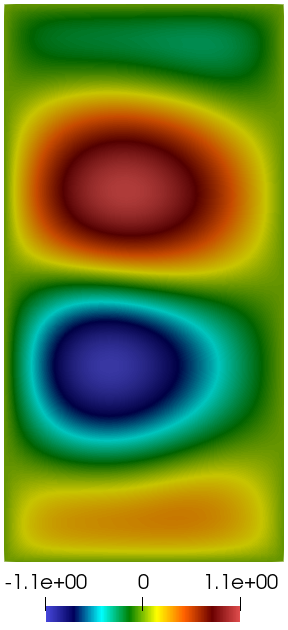} & \hspace{-0.4cm}\includegraphics[align=c,scale = 0.25]{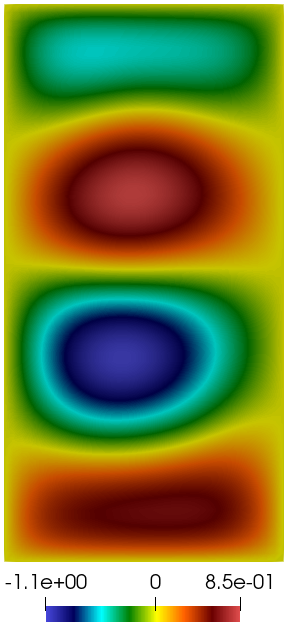} & \hspace{-0.4cm}\includegraphics[align=c,scale = 0.25]{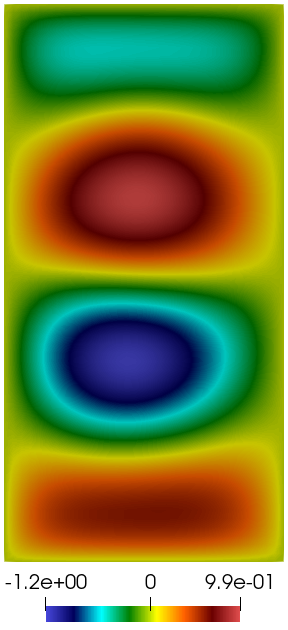} \\
        \hspace{-0.4cm}$\tilde{q}$  & \hspace{-0.4cm}\includegraphics[align=c,scale = 0.25]{Exp2_figs/QGE/qMean_DNS.png} & \hspace{-0.4cm}\includegraphics[align=c,scale = 0.25]{Exp2_figs/QGE-noSGS/qMean_qge32.png} & \hspace{-0.4cm}\includegraphics[align=c,scale = 0.25]{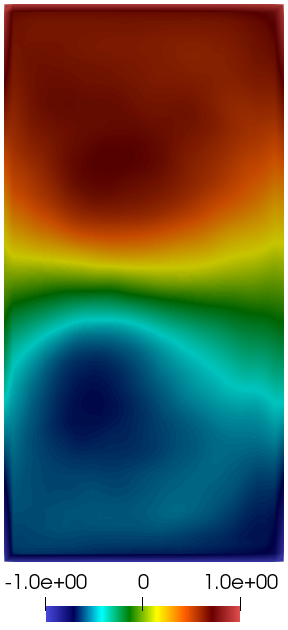} & \hspace{-0.4cm}\includegraphics[align=c,scale = 0.25]{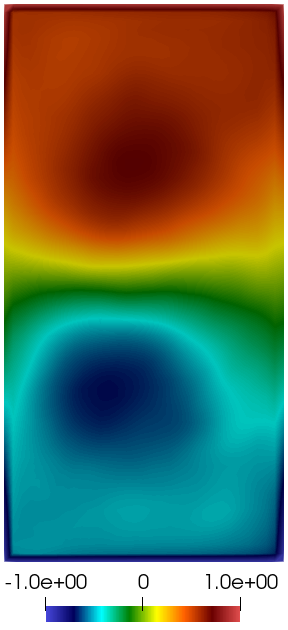} & \hspace{-0.4cm}\includegraphics[align=c,scale = 0.25]{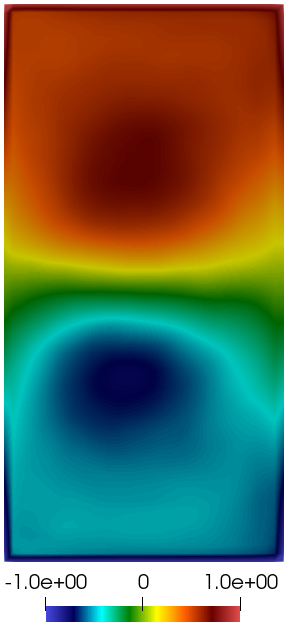}
    \end{tabular}
    \caption{QGE-EFR, experiment 2: Time-averaged stream function $\tilde{\psi}$ (first row) and potential vorticity $\tilde{q}$ (second row) computed from the DNS (first column), QGE with no SGS model (second column), QGE-EFR model with indicator functions $a_{lin}$ (third column), $a_{D}$ (fourth column), and $a_S$ (fifth column). The coarse mesh size is $h = 1/32$ and we set $\alpha = 2h$.}
    \label{fig:Mean_32_2h_exp2_efr}
\end{figure}

Next, we consider experiment 2. 
Fig.~\ref{fig:Mean_32_2h_exp2_efr}
compares the time averaged stream function $\tilde{\psi}$ and potential vorticity $\tilde{q}$ computed by the DNS (the first column), QGE without filtering and QGE-EFR models.
Again, the corse mesh has size $h = 1/32$, but, given challenge of a higher Reynolds number in
experiment 2, the setting used so far of $\alpha = h$ provides insufficient artificial dissipation
and it is thus changed to $\alpha = 2h$.
In fact, only the linear indicator function $a_{lin}$ introduces enough dissipation for the computed stream function $\tilde{\psi}$ to exhibit the correct four-gyre structure when  $\alpha = h$. With $\alpha = 2h$, the QGE-EFR solutions are in agreement agreement with the DNS, correctly reproducing the four-gyre structure in $\tilde{\psi}$ and the large-scale features of $\tilde{q}$ with no spurious oscillations at the boundary. We also observe that among the QGE-EFR solutions, the ones computed with $a_D$ and $a_S$ produce results in better qualitative agreement with the DNS.

\begin{figure}[htb!]
    \centering
    \begin{subfigure}{0.45\linewidth}
        \includegraphics[width = \linewidth]{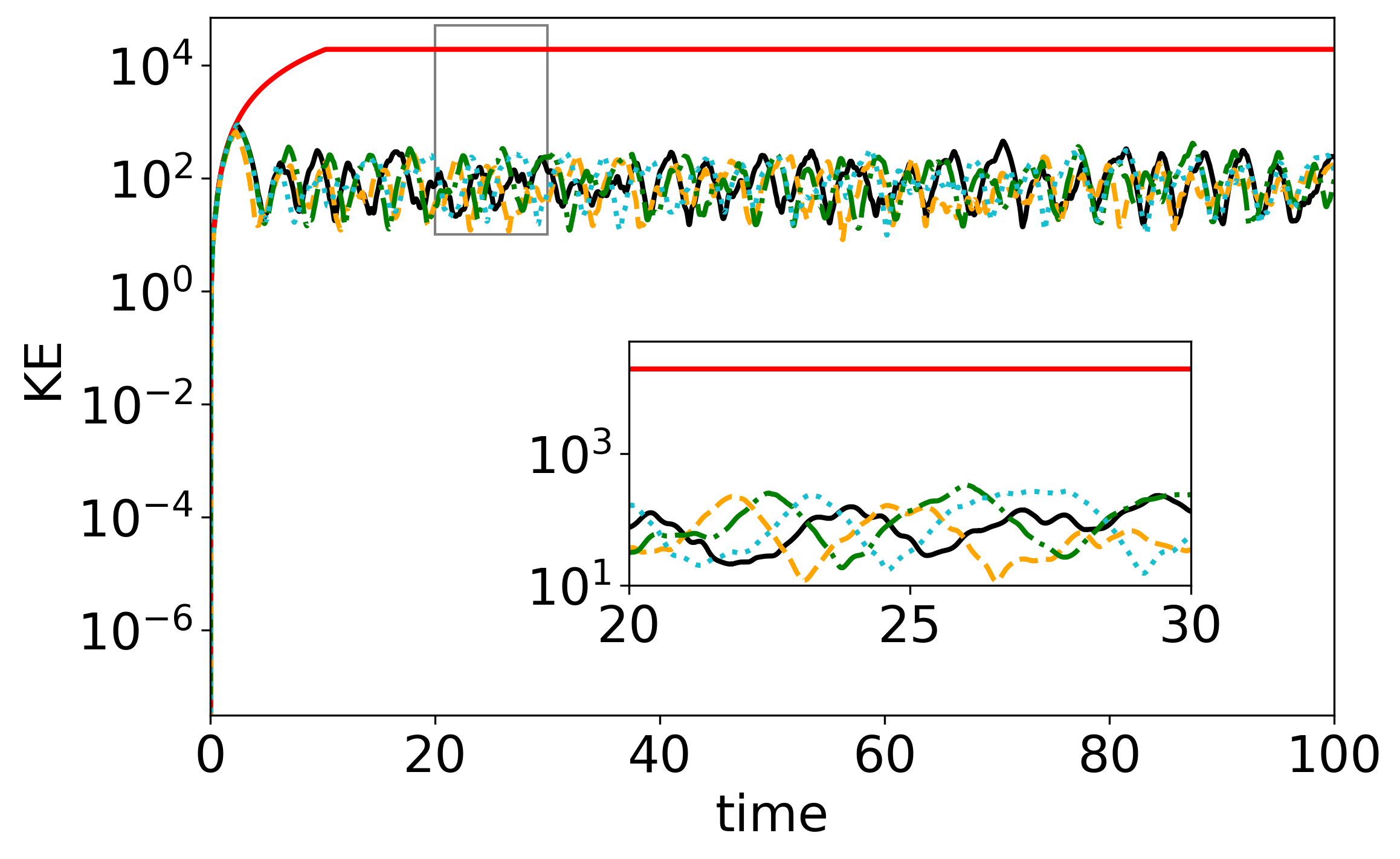} 
        \caption{Kinetic energy}
    \end{subfigure}
    \begin{subfigure}{0.45\linewidth}
        \includegraphics[width = \linewidth]{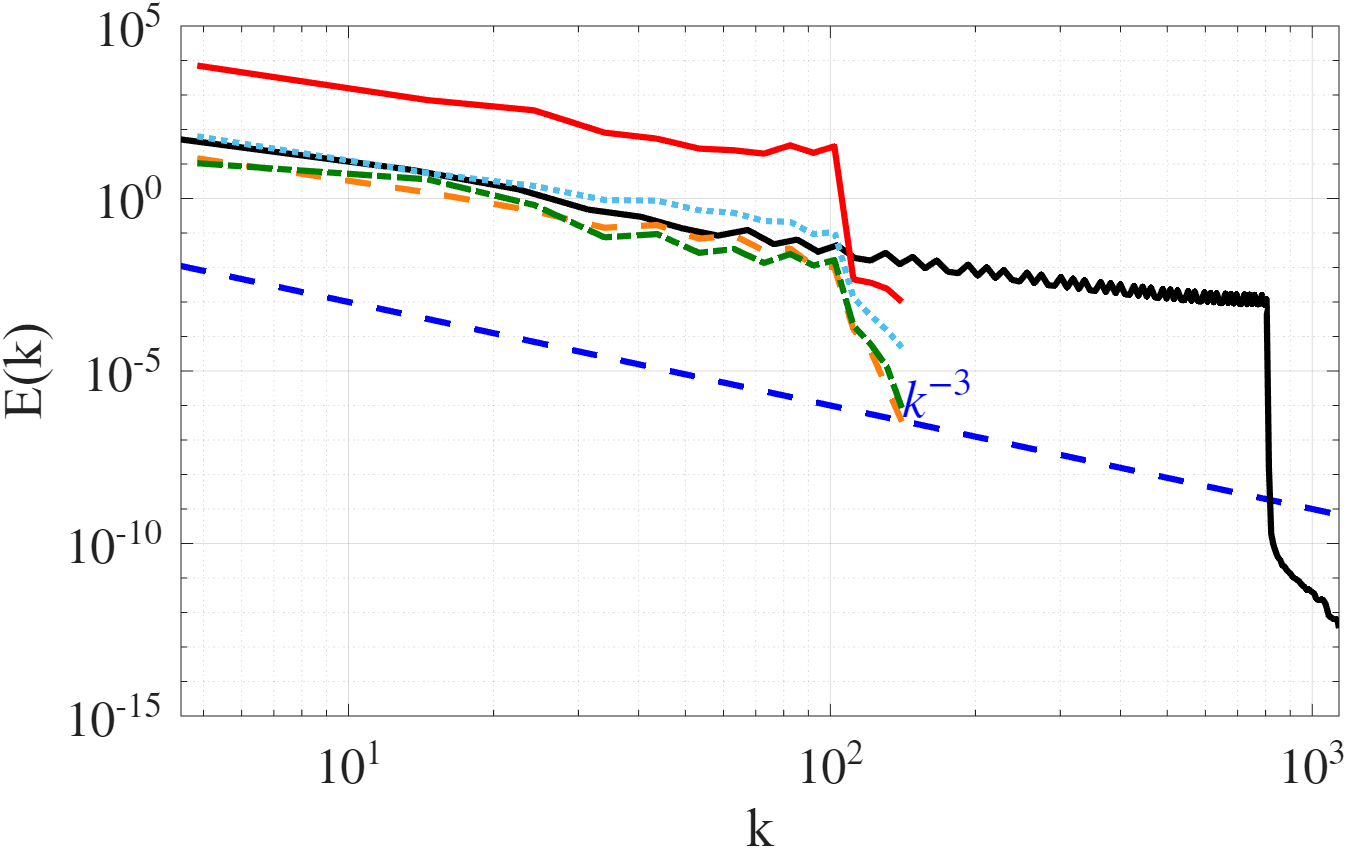}
        \caption{KE spectra}
    \end{subfigure}
    \begin{subfigure}{0.7\linewidth}
        \includegraphics[width = \linewidth]{Exp1_figs/EFR/legend.png}
    \end{subfigure}
    \caption{QGE-EFR, experiment 2: (a) Evolution of the KE and (b) KE spectrum at time $t = 20$ given by the DNS, the
    QGE with no SGS model,  and the QGE-EFR model with indicator functions $a_{lin}$, $a_D$, and $a_S$ for mesh $h = 1/32$. We set $\alpha = 2h$.}
    \label{fig:ke_kespec_32_exp2_efr}
\end{figure}

Fig.~\ref{fig:ke_kespec_32_exp2_efr} reports the evolution of the kinetic energy and the corresponding kinetic energy spectra for experiment 2 using the coarse mesh $h=1/32$ and filtering radius $\alpha=2h$. Consistent with the observations from Fig.~\ref{fig:Mean_32_2h_exp2_efr}, the coarse-mesh QGE simulation without filtering overpredicts the kinetic energy given by the DNS. In contrast, all QGE-EFR simulations provide a stable energy evolution and reproduce the mean KE level with good accuracy. Among the tested indicators, the deconvolution-based indicator $a_D$ yields the closest agreement with the DNS throughout the simulation, while the linear indicator $a_{lin}$, although more dissipative, remains effective in controlling the energy growth. The energy spectra shown in Fig.~\ref{fig:ke_kespec_32_exp2_efr} (b) further confirm these trends: all QGE-EFR models recover the correct distribution of energy over the resolved scales and significantly improve upon the unfiltered QGE solution. %In particular, $a_D$ provides the best match to the DNS spectrum at intermediate and high wavenumbers, indicating a more accurate representation of the energy transfer mechanisms. 
Overall, the results demonstrate that increasing the filtering radius to $\alpha=2h$ supplies sufficient stabilization for the challenging high-Reynolds number regime of experiment 2 while preserving the main flow statistics.

The numerical results in this subsection suggest that the QGE-EFR model is accurate provided that the mesh is not too coarse and the filter introduces sufficient dissipation. Among the indicator functions we considered, $a_D$ provides the most robust performance. For the more challenging experiment 2, $a_{lin}$ introduces additional dissipation that may be beneficial, making it the only indicator capable of recovering the correct four-gyre structure for the case with mesh size $h = 1/32$ and $\alpha = h$. These results suggest that $a_D$ is the preferred choice for EFR in general, whereas $a_{lin}$ may be advantageous for highly under-resolved flows requiring additional stabilization.

\subsection{A more realistic test case} \label{sec:mediterranean}

\begin{figure}[htb!]
    \centering
        \includegraphics[width = .8\linewidth]{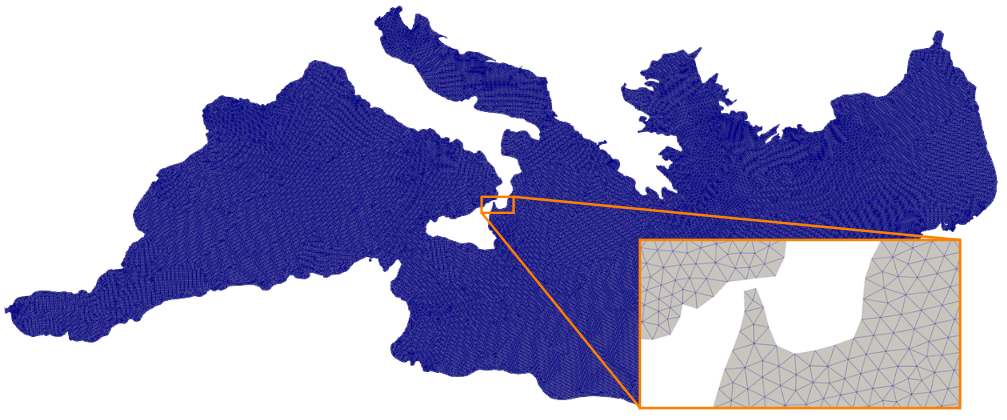}
        \caption{Fine mesh of the Mediterranean Sea.} %\sachin{ How does this look? If this doesn't look good, we can make a separate picture for the fine mesh as Prof Quaini suggested.}}
    \label{fig:med_mesh_fine}
\end{figure}

\begin{figure}[htb!]
\centering
        \includegraphics[width = .5\linewidth]{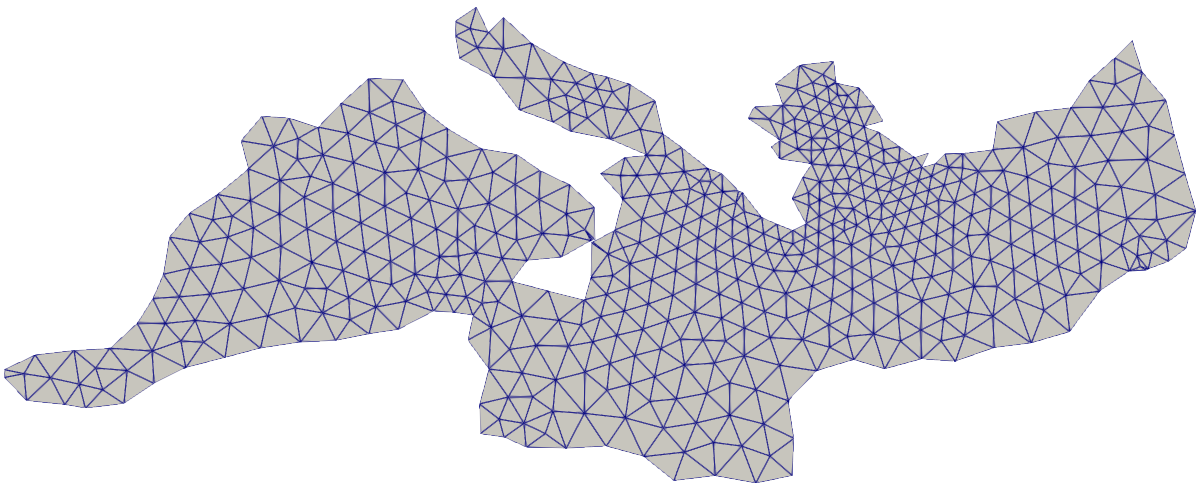}
        \caption{Coarse mesh of the Mediterranean Sea.}
    \label{fig:med_mesh_coarse}
\end{figure}

To assess the performance of the BV-Bardina and QGE-EFR models on a realistic benchmark, we consider the geometry of the Mediterranean Sea. 
The coastline, obtained from
the  Global Self-consistent, Hierarchical, High-resolution Geography (GSHHS) database \cite{GSHHS}, was kindly made available to us by the authors of \cite{QGE-review,Foster2013,foster2013finite}. 
% \anna{We should cite Traian's papers where they use it and say it was kindly made available to us by the authors}. \sachin{It has been used in \cite{QGE-review}, they have shown  the mesh but did not run any simulation on it, the authors say that it has been used in \cite{Foster2013} and in a related thesis \cite{foster2013finite} but in  \cite{Foster2013} and \cite{foster2013finite}, they  have used the simplified geometry not the one  we have used.} 
In contrast to the benchmarks in the previous subsections, this test involves a computational domain with complex features. The computational mesh is generated with GMSH \cite{GMSH}. 
Several major islands have been removed and Sicily has been attached to mainland Italy to obtain a monoconnected domain, which ensures
a unique stream function. See \cite{myers1995diagnostic,gunzburger2012finite,gunzburger1988finite,gunzburger1988finite2} for a discussions on this. 
We consider a fine mesh with 
size $h = 1.22e-03$ (leading to
96,381 cells) shown in  Fig.~\ref{fig:med_mesh_fine}, 
and a coarse mesh with size 
$h = 1.56e-02$ (leading to 986 cells) reported in 
Fig.~\ref{fig:med_mesh_coarse}.
Considering a characteristic length of $L = 1000$ Km \cite{Foster2013}, the dimensional mesh sizes are 
$1.22$ Km and $15.6$ Km, respectively. 

While our previous benchmarks in Sec.~\ref{sec:bv-bardina}-\ref{sec:qge-efr} employ the QGE in the potential vorticity-stream function ($q$-$\psi$) formulation presented in Sec.~\ref{sec:QGE}, for this test we use the equivalent vorticity-stream function ($\omega$-$\psi$) formulation. See, e.g., \cite{Monteiro2015}. 
This choice is so that we can rewrite boundary condition
\eqref{eq:BV5_comp2} as $\omega = 0$, which allows us to avoid finding the center of the basin and non-dimensionalizing the $y$-coordinate. The BV-Bardina and
QGE-EFR models have a straightforward adaptation to the ($\omega$-$\psi$) formulation.

Most of the settings for this test are taken from \cite{Foster2013}, including the wind forcing prescribed as $F = \sin\left(\pi y/4\right)$. %The eddy viscosity $A$ is taken to be $2.0\times 10^3$ m$^2$~s$^{-1}$ and the gradient of the Coriolis frequency as $\beta = 1.742\times10^{-11}$ m$^{-1}$~s$^{-1}$. 
In contrast to \cite{Foster2013}, we take the following values of $Ro$ and $Re$:
\begin{equation*}
    Ro = 6.74\times 10^{-5}, \quad \text Re = 527,
\end{equation*}
which are significantly more challenging than the ones used in \cite{Foster2013}. %This means that the characteristic velocity is $U = \frac{Re~A}{L} = 1.054$ m~s$^{-1}$. 
The nondimensional time step is taken as $\Delta t = 2.5e-05$, which corresponds to approximately $24$ s in dimensional time. 
The nondimensional time interval for the simulation is
$(0, 5]$ and 
%Consequently, the time interval of interest $[1, 5]$  %\lander{more accurately, 43.89} %\sachin{we are taking time averages in the interval [1,5]} 
%Throughout this section, 
the time-averaged stream function $\tilde{\psi}$ and vorticity $\tilde{\omega}$ are computed over the time interval $[1,5]$, corresponding to about $44$ days. Lastly, we note that the Munk scale for this problem is $\delta_M/L = 5.03e-03$ and thus, the fine mesh size is almost 5 times smaller than the Munk scale, while the coarse mesh size is about three times larger.

\begin{figure}[htb!]
    \centering
    \begin{subfigure}{0.45\linewidth}
        \centering
        \includegraphics[width = \linewidth]{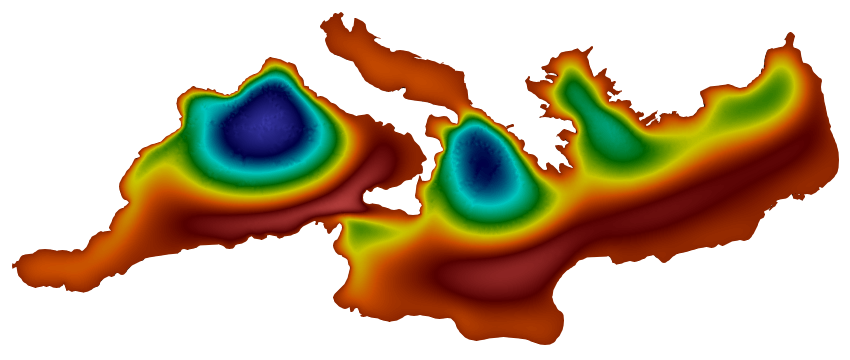}
        \caption{DNS}
    \end{subfigure}
    \begin{subfigure}{0.45\linewidth}
        \centering
        \includegraphics[width = \linewidth]{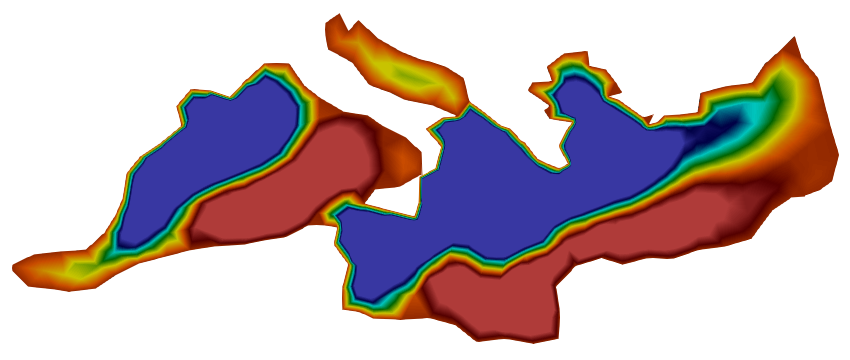}
        \caption{BV-Bardina ($a_{lin}$)}
    \end{subfigure}
    \begin{subfigure}{0.45\linewidth}
        \centering
        \includegraphics[width = \linewidth]{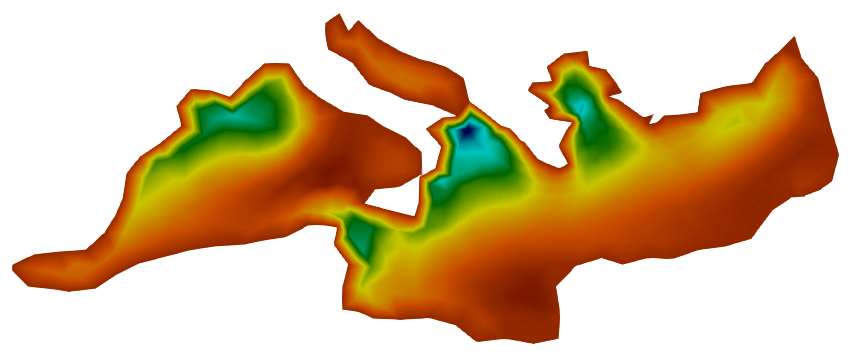}
        \caption{QGE-EFR ($a_{lin}$)}
    \end{subfigure}
    \begin{subfigure}{0.45\linewidth}
        \centering
        \includegraphics[width = \linewidth]{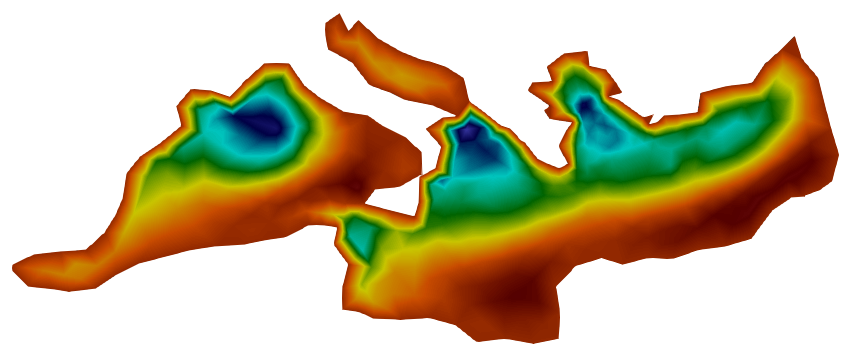}
        \caption{QGE-EFR ($a_D$)}
    \end{subfigure}
    \begin{subfigure}{0.45\linewidth}
        \centering
        \includegraphics[width = \linewidth]{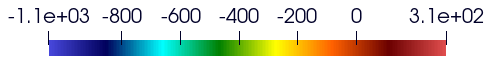}
    \end{subfigure}
    \caption{Mediterranean Sea: Time-averaged vorticity  $\tilde{\omega}$ computed by the (a) DNS, (b) BV-Bardina model with $a_{lin}$, (c) QGE-EFR model with $a_{lin}$, and (d) QGE-EFR model with $a_D$. All simulations with regularization are performed with $\alpha = 1.7h$.}
    \label{fig:omega_ave}
\end{figure}

\begin{figure}[htb!]
    \centering
    \begin{subfigure}{0.45\linewidth}
        \centering
        \includegraphics[width = \linewidth]{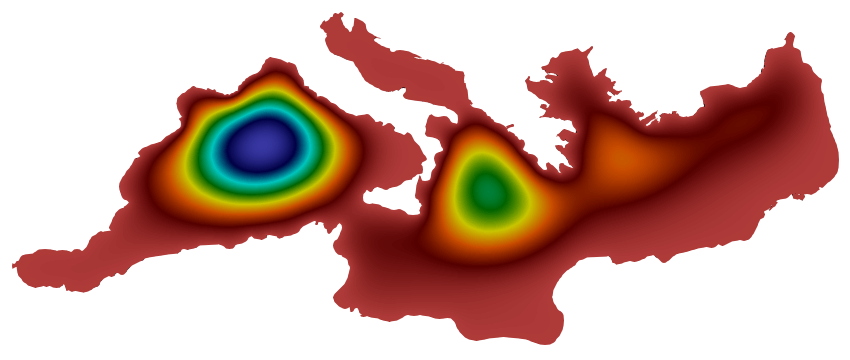}
        \caption{DNS}
    \end{subfigure}
    \begin{subfigure}{0.45\linewidth}
        \centering
        \includegraphics[width = \linewidth]{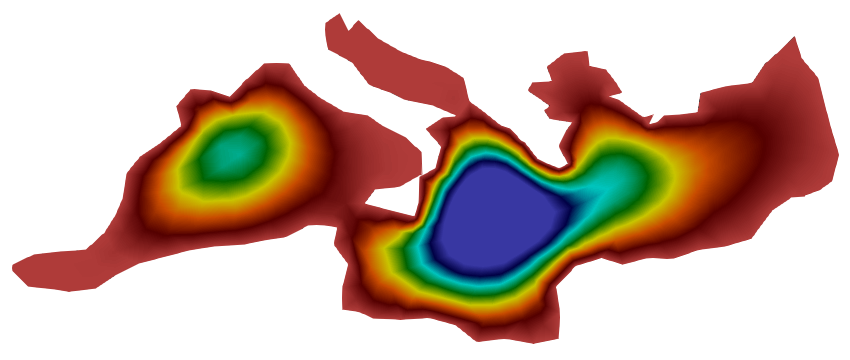}
        \caption{BV-Bardina ($a_{lin}$)}
    \end{subfigure}
    \begin{subfigure}{0.45\linewidth}
        \centering
        \includegraphics[width = \linewidth]{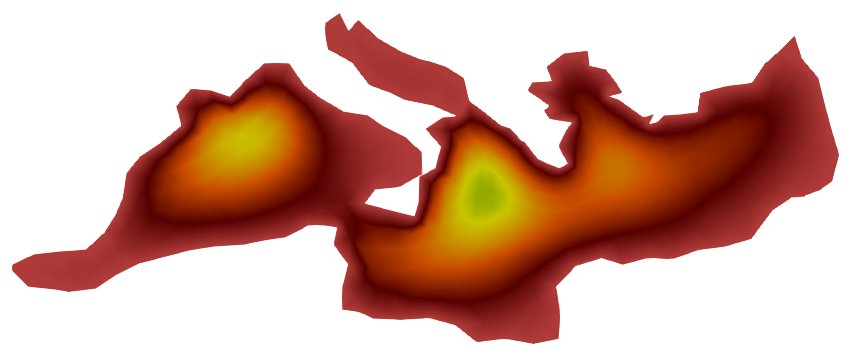}
        \caption{QGE-EFR ($a_{lin}$)}
    \end{subfigure}
    \begin{subfigure}{0.45\linewidth}
        \centering
        \includegraphics[width = \linewidth]{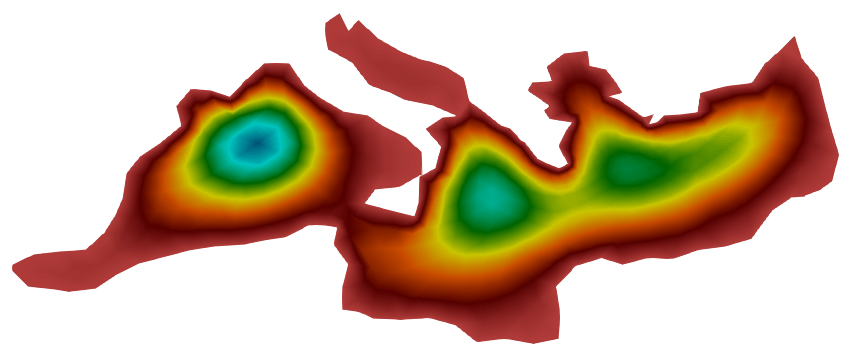}
        \caption{QGE-EFR ($a_D$)}
    \end{subfigure}
    \begin{subfigure}{0.45\linewidth}
        \centering
        \includegraphics[width = \linewidth]{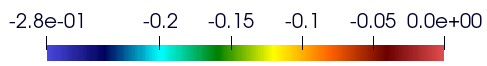}
    \end{subfigure}
    \caption{Mediterranean Sea: Time-averaged stream function  $\tilde{\psi}$ computed by the (a) DNS, (b) BV-Bardina model with $a_{lin}$, (c) QGE-EFR model with $a_{lin}$, and (d) QGE-EFR model with $a_D$. All simulations with regularization are performed with $\alpha = 1.7h$.}
    \label{fig:psi_ave}
\end{figure}

We report $\tilde{\omega}$ and $\tilde{\psi}$ given by the DNS in Fig.~\ref{fig:omega_ave} (a) and Fig.~\ref{fig:psi_ave} (a), respectively. In Fig.~\ref{fig:omega_ave} (a), we see 
two large regions with negative
vorticity (in dark blue), the larger being West of Italy and the smaller between Southern Italy and Western Greece.
In the same regions, we observe 
two large gyres in 
Fig.~\ref{fig:psi_ave} (a).
For the coarse mesh simulations, we fix the filtering radius to $\alpha = 1.7h$. We note that the coarse simulation of the QGE with no SGS and BV-Bardina model with $a_S$ become unstable in the time interval of interest. %Consequently, the time-averaged vorticity $\tilde{\omega}$ and streamfunction $\tilde{\psi}$ are not included in the comparison, and only stable simulations are presented below.
Among the stable simulations, the BV-Bardina model with $a_{lin}$ fails to reproduce the shape and size of the dominant gyres observed in the DNS. Compare Fig.~\ref{fig:psi_ave} (b) with Fig.~\ref{fig:psi_ave} (a). In addition, 
the vorticity field $\tilde{\omega}$ in Fig.~\ref{fig:omega_ave} (b) presents regions with large positive (dark red) and large negative (dark blue) values, suggesting that the model does not effectively dissipate unphysical oscillations. 
In contrast, the QGE-EFR model with $a_{lin}$ remains stable and the shape and size of the dominant gyres are qualitatively preserved (compare Fig.~\ref{fig:psi_ave} (c) with Fig.~\ref{fig:psi_ave} (a)), 
but it exhibits excessive artificial dissipation.  Replacing the linear indicator with the deconvolution-based indicator $a_D$ significantly improves the solution obtained with the QGE-EFR model. The two dominant gyres are recovered at approximately the correct locations and with magnitudes closer to those of the DNS. See Fig.~\ref{fig:omega_ave} (d) and \ref{fig:psi_ave} (d).

\begin{figure}[htb!]
    \centering
    \begin{subfigure}{0.45\linewidth}
        \centering
        \includegraphics[width = \linewidth]{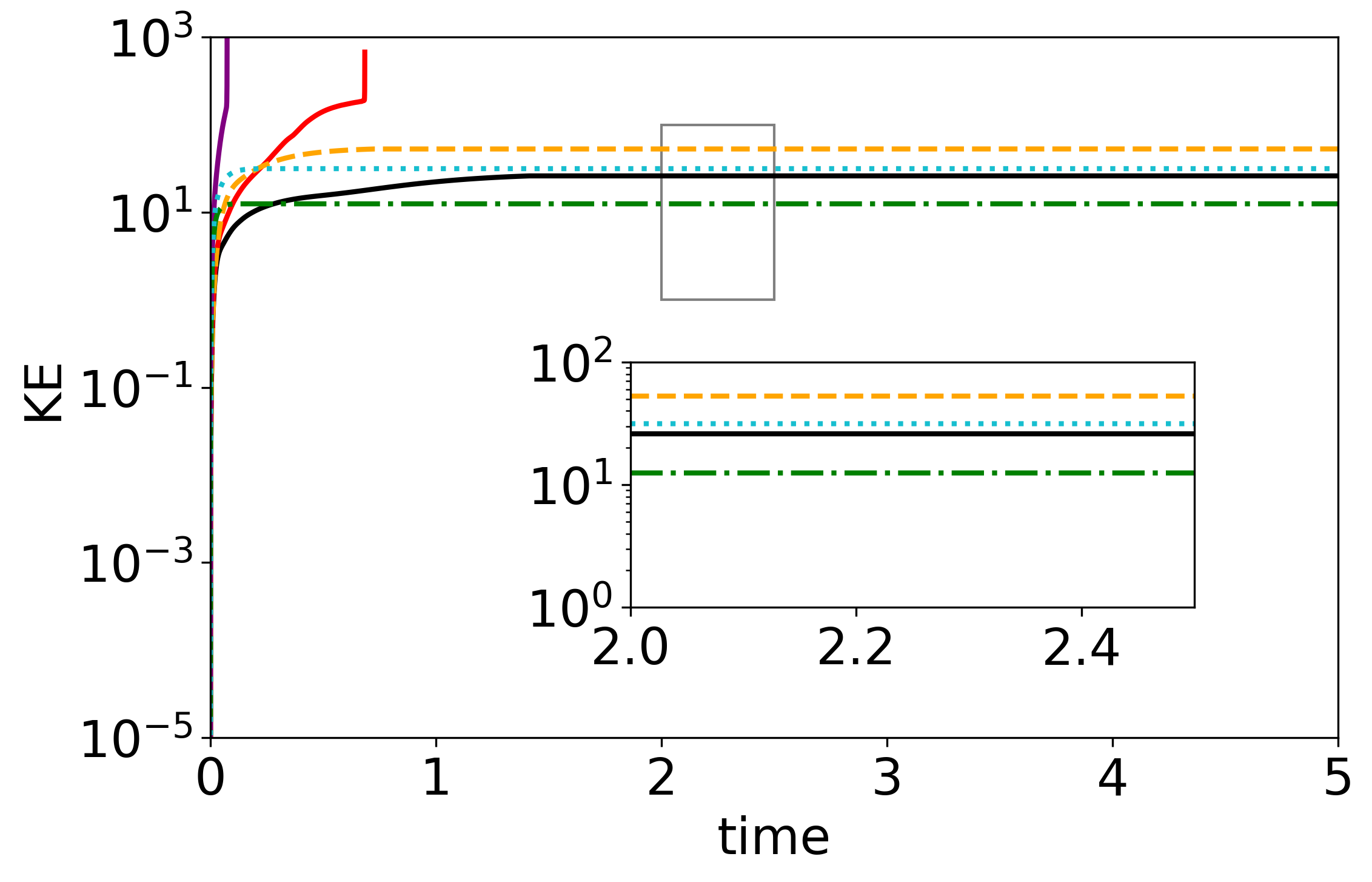}
        \caption{Kinetic energy}
    \end{subfigure}
    \begin{subfigure}{0.45\linewidth}
        \centering
        \includegraphics[width = \linewidth]{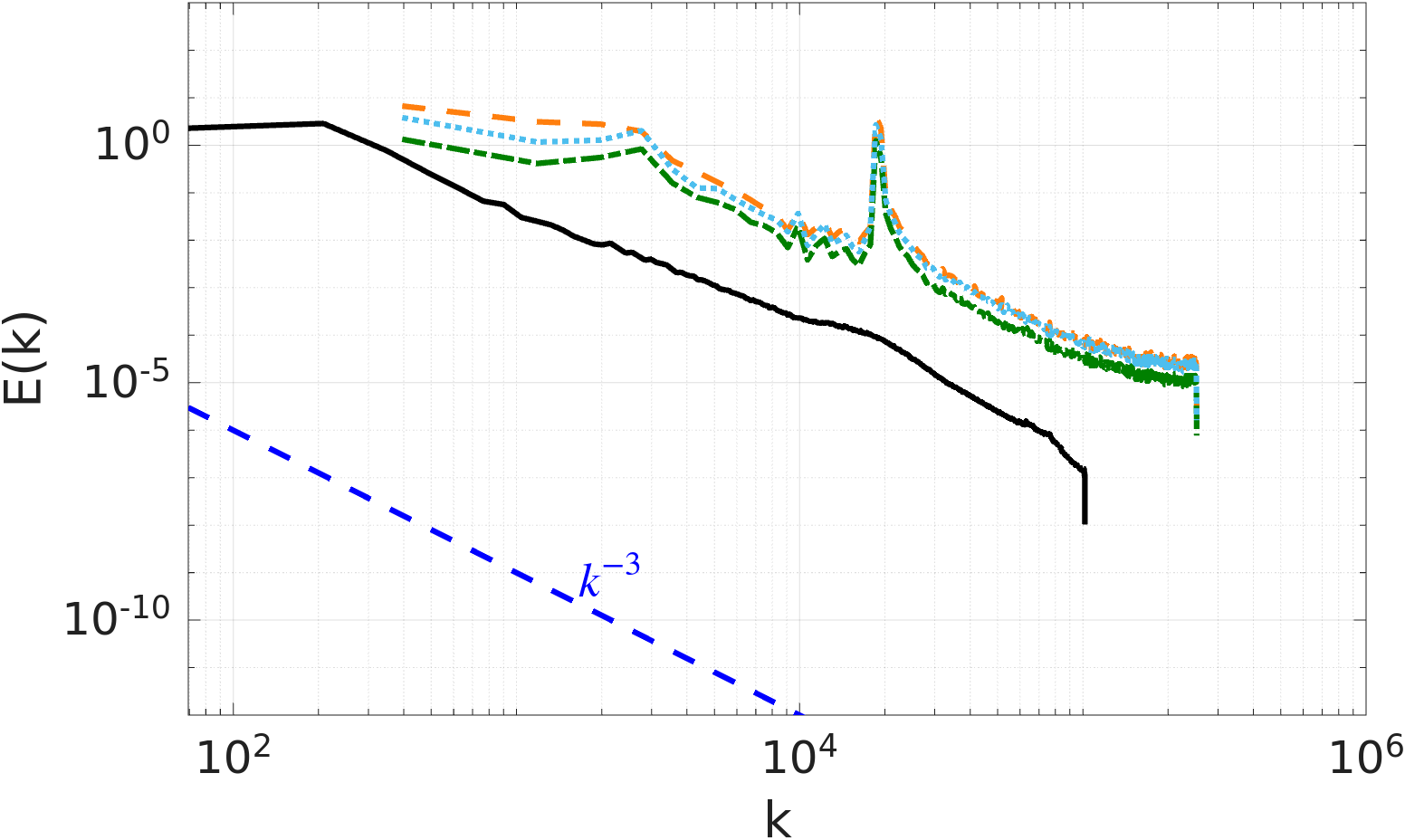}
        \caption{KE spectra}
    \end{subfigure}
    \begin{subfigure}{0.85\linewidth}
        \centering
        \includegraphics[width = \linewidth]{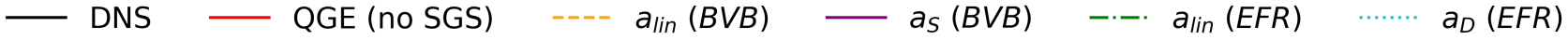}
    \end{subfigure}
    \caption{Mediterranean Sea: (a) Evolution of the KE and (b) KE spectrum at time $t = 2$ given by the DNS, the
    QGE with no SGS model, BV-Bardina (BVB) with $a_{lin}$ and $a_{S}$, and QGE-EFR model with $a_{lin}$ and $a_{D}$. We set $\alpha = 1.7h$.}
    \label{fig:energy}
\end{figure}

Fig.~\ref{fig:energy} compares the evolution of the kinetic energy and the corresponding KE spectrum at $t = 2$. In Fig.~\ref{fig:energy} (a), we report also the KE evolution for the coarse QGE simulation with no SGS and BV-Bardina model with $a_S$ to show the numerical instability. All other models in Fig.~\ref{fig:energy} (a) remain stable throughout the time interval of interest. The value of KE given by the
BV-Bardina model with $a_{lin}$ and the QGE-EFR model with $a_D$ is larger than the value given by the DNS, while the value given by the QGE-EFR model with $a_{lin}$ is smaller, confirming the more dissipative behavior already seen in Fig.~\ref{fig:omega_ave} (c). Among these three models, the QGE-EFR model with $a_D$ produces the KE evolution in closest agreement with the DNS. As for the spectra in Fig.~\ref{fig:energy}, the largest discrepancies between the coarse-grid models and the DNS occur at the highest resolved wavenumbers, where under-resolution effects are expected to be strongest. Nevertheless, the dominant large-scale energy-containing structures are well captured by the QGE-EFR model. 
%All stable simulations produce KE spectra which have similar decay as the DNS with respect to increasing wavenumber $k$. However, the coarse simulations with SGS models retain more energy over a broad range of wavenumbers than the DNS, particularly at intermediate and high wavenumbers. Among the methods considered, the two QGE-EFR models provide the closest agreement with the DNS spectrum. By comparison, the BV-Bardina model consistently predicts higher energy levels across the resolved scales, suggesting that BV-Bardina does not introduce sufficient artificial dissipation to the system.

\begin{figure}[htb!]
    \centering
    \begin{subfigure}{0.45\linewidth}
        \centering
        \includegraphics[width = \linewidth]{Mediterranean/DNS/psiMean.png}
        \caption{DNS}
    \end{subfigure}
    \begin{subfigure}{0.45\linewidth}
        \centering
        \includegraphics[width = \linewidth]{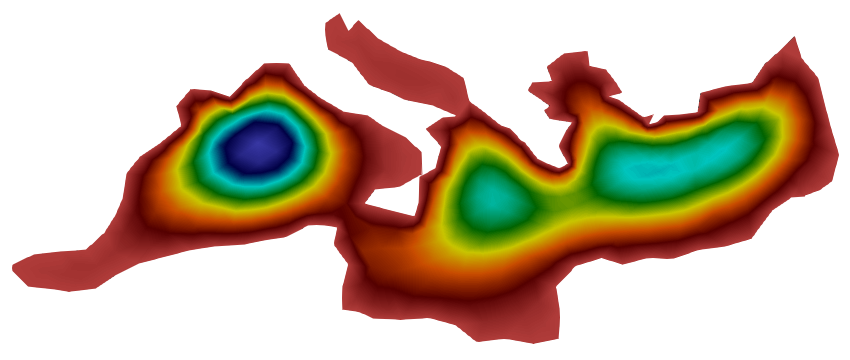}
        \caption{QGE-EFR ($a_D$), $\alpha = \sqrt{2}h$}
    \end{subfigure}
    \begin{subfigure}{0.45\linewidth}
        \centering
        \includegraphics[width = \linewidth]{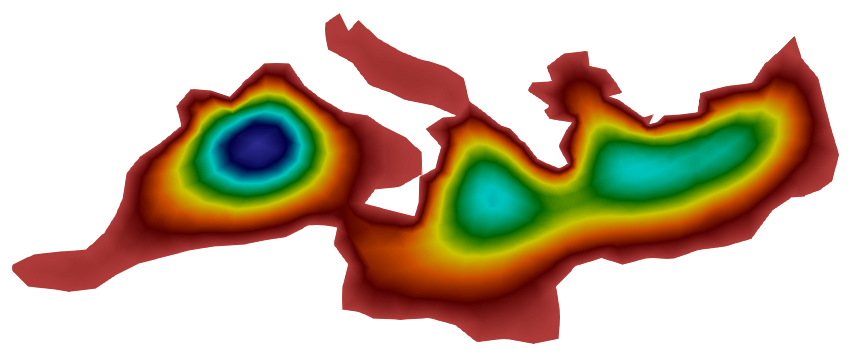}
        \caption{QGE-EFR ($a_D$), $\alpha = 1.5h$}
    \end{subfigure}
    \begin{subfigure}{0.45\linewidth}
        \centering
        \includegraphics[width = \linewidth]{Mediterranean/EFR_deconvu/alpha_1_7h/psiMean.png}
        \caption{QGE-EFR ($a_D$), $\alpha = 1.7h$}
    \end{subfigure}
    \begin{subfigure}{0.45\linewidth}
        \centering
        \includegraphics[width = \linewidth]{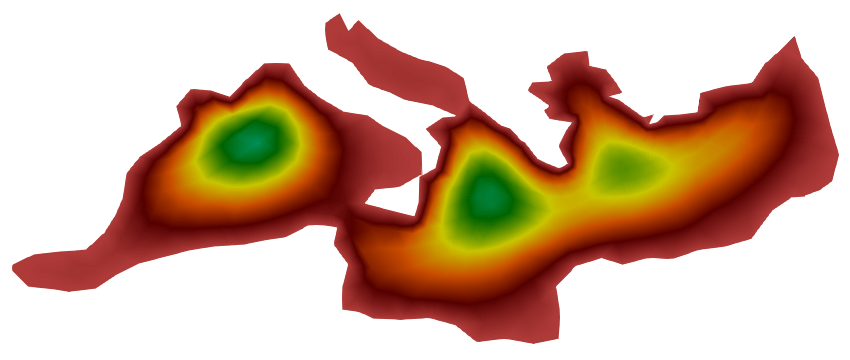}
        \caption{QGE-EFR ($a_D$), $\alpha = 2h$}
    \end{subfigure}\\
    \begin{subfigure}{0.45\linewidth}
        \centering
        \includegraphics[width = \linewidth]{Mediterranean/DNS/legend_psiMean.png}
    \end{subfigure}
    \caption{Mediterranean Sea: Time-averaged stream function $\tilde{\psi}$ computed by the (a) DNS and (b)-(e) QGE-EFR model with $a_D$ and different values of $\alpha$.}
    \label{fig:psi_ave_efr_sens}
\end{figure}

Since the QGE-EFR model with $a_D$ provides the closest overall agreement with the DNS, next
we investigate the sensitivity of this model to the filtering radius $\alpha$. In particular, we examine how varying $\alpha$ influences the large-scale structures. Fig.~\ref{fig:psi_ave_efr_sens} displays the time-averaged stream function $\tilde{\psi}$ given by the QGE-EFR model with $a_D$ and $\alpha = \sqrt{2}h, 1.5h, 1.7h$, and $2h$.
Since increasing $\alpha$ corresponds to increase in artificial dissipation, the maximum and minimum values of $\tilde{\psi}$ gradually reduce as $\alpha$ is increased. While the solutions obtained with $\alpha = \sqrt{2}h, 1.5h$ are in close agreement with the DNS, $\alpha = 2h$ produces excessive smoothing. %Thus, these results suggest that QGE-EFR exhibits only moderate sensitivity to the filtering radius $\alpha$. 
We also notice that, although
the dominant gyres are preserved over the range of values considered, the QGE-EFR model does show a certain sensitivity to the choice of the filtering radius.

%\anna{Please also add a picture of the fine mesh.}\sachin{done.}

%\subsection{Computational cost comparison} \label{sec:comp_cost}

Finally, we report on the computational costs for this challenging test. All simulations are performed on a computing server with Intel\textsuperscript{\textregistered} Xeon\textsuperscript{\textregistered}, 256 Gb of RAM, and a 64-bit version of Linux. The DNS %, performed on a computational mesh with 96,381 cells,
takes about 2 h 18 min to complete. %We recall that the coarse simulation of the QGE with no SGS and BV-Bardina with $a_S$ became unstable. 
Both the BV-Bardina and QGE-EFR models with $a_{lin}$ take about 1.5 min. This corresponds to a speedup of about 85 times compared to the DNS. On the other hand, the QGE-EFR model with the indicator function $a_D$ requires 2 min 9 s, which is about 64 times faster than the DNS. Although the QGE-EFR model with $a_D$ is slightly more expensive than the linear models, it provides
the closest results to the DNS at a fraction of its computational cost.

\begin{comment}
\begin{table}[h!]
    \centering
    \begin{tabular}{|lcc|}
        \hline
        Model & CPU Time & Speed up factor \\
        \hline
        DNS & 2h 18m & - \\
        QGE with no SGS & - & - \\
        BV-Bardina with $a_{lin}$ & 1m 37s & 85 \\
        BV-Bardina with $a_S$ & - & - \\
        QGE-EFR with $a_{lin}$ & 1m 38s & 85 \\
        QGE-EFR with $a_{D}$ & 2m 9s & 64 \\
        \hline
    \end{tabular}
    \caption{Caption}
    \label{tab:placeholder}
\end{table}
\end{comment}

\section{Conclusion}\label{sec:conclusion}
We introduced the Evolve-Filter-Relax framework for the under-resolved quasi-geostrophic equations and argued that it can be interpreted as a splitting scheme for a perturbed QGE system with additional dissipative terms.
We showed that the EFR algorithm is equivalent to an eddy viscosity model in LES and we considered three indicator functions to tune the amount and location of eddy viscosity: a constant indicator function, an indicator functions that recover a Smagorinsky-like model, and an indicator function based on an approximate deconvolution operator.
In addition, a practical expression for the relaxation parameter was derived based on physical and discretization parameters.
The performance of the EFR algorithm for the QGE
was assessed against a nonlinear variant of the BV-Bardina model, which also uses different indicator functions. 

Numerical results on the double-gyre wind forcing benchmark showed that both the 
QGE-EFR and the nonlinear BV-Bardina model
significantly improve the accuracy and stability of the of coarse mesh simulations compared with the QGE with no sub-grid scale model.
Indeed, they successfully recover time-averaged flow structures, realistic kinetic energy levels, and accurate energy spectra, even when using meshes much coarser than those required by a DNS.
Among the investigated models, the QGE-EFR model with the deconvolution-based indicator delivered the best balance between accuracy, stability, and computational efficiency.
In addition, it
performed well in the realistic Mediterranean Sea test case, where it accurately captured 
both the large-scale gyres and the evolution
of the kinetic energy, while being 
significantly less expensive (about 64 times)
than a DNS.

\section*{Acknowledgement}
The authors sincerely thank Traian Iliescu for providing the computational mesh of the Mediterranean Sea, and Michele Girfoglio for his assistance in computing and visualizing the kinetic energy spectra and for insightful discussions that helped improve this work. 

As we were completing this manuscript, our mentor, collaborator, and friend Alessandro Veneziani passed away unexpectedly and far too soon. This work is a direct continuation of the ideas in \cite{BQV}, shaped in no small part by his insight and vision. We gratefully acknowledge his lasting influence on our work and careers and dedicate this paper to his memory.

\section*{Data statement}
The data that support the findings of this study are available from the corresponding author upon reasonable request.

\bibliography{QGE} 

%\appendix
%\section{The EFR algorithm as a solver for the QGE with perturbations} \label{sec:append_perturb}

\end{document}